\documentclass[a4paper, 11pt]{article}
\usepackage{amsmath}
\usepackage{amsfonts}
\usepackage{amssymb}
\usepackage[english]{babel}
\usepackage{graphicx}
\usepackage{amsthm}
\usepackage[title]{appendix}

\usepackage{longtable}
\usepackage{tabularx}
\allowdisplaybreaks

\newtheorem{theorem}{Theorem}[section]

\newtheorem{lemma}[theorem]{Lemma}

\newtheorem{construction}[theorem]{Construction}
\def\whitebox{{\hbox{\hskip 1pt
 \vrule height 6pt depth 1.5pt
 \lower 1.5pt\vbox to 7.5pt{\hrule width
    3.2pt\vfill\hrule width 3.2pt}%
 \vrule height 6pt depth 1.5pt
 \hskip 1pt } }}
\def\qed{\ifhmode\allowbreak\else\nobreak\fi\hfill\quad\nobreak
     \whitebox\medbreak}

\newcommand{\ignore}[1]{}

\begin{document}
\baselineskip 16pt
\title{The Hamilton-Waterloo Problem for $C_3$-factors and $C_n$-factors}

\author{\small   Li Wang, Fen Chen and Haitao Cao \thanks{Research
supported by the National Natural Science Foundation of China
under Grant 11571179, the Natural Science Foundation of Jiangsu
Province under Grant No. BK20131393, and the Priority Academic
Program Development of Jiangsu Higher Education Institutions.
E-mail: {\sf caohaitao@njnu.edu.cn.}} \\
\small Institute of Mathematics, \\ \small  Nanjing Normal
University, Nanjing 210023, China}

\date{}
\maketitle
\begin{abstract}
The Hamilton-Waterloo problem asks for a 2-factorization of $K_v$
(for $v$ odd) or $K_v$ minus a $1$-factor (for $v$ even) into
$C_m$-factors and $C_n$-factors. We completely solve the
Hamilton-Waterloo problem in the case of $C_3$-factors and
$C_n$-factors for $n=4,5,7$.

\smallskip
\end{abstract}

\section{Introduction}

In this paper, the vertex set and the edge set of a graph $H$ will
be denoted by $V(H)$ and $E(H)$, respectively.   We denote the
cycle of length $k$ by $C_k$, the complete graph on $v$ vertices
by $K_v$, and the complete $u$-partite graph with $u$ parts of
size $g$ by $K_u[g]$. A {\it factor} of a graph $H$ is a spanning
subgraph of $H$. Suppose $G$ is a subgraph of a graph $H$, a {\it
$G$-factor} of $H$ is a set of edge-disjoint subgraphs of $H$,
each isomorphic to $G$. And a {\it $G$-factorization} of $H$ is
a set of edge-disjoint $G$-factors of $H$.
Many authors \cite{AH,ASSW, HS, JLA, L, LD,PWL, R} have contributed to
prove the following result.

\begin{theorem}\label{Kv}
There exists a $C_k$-factorization of $K_u[g]$ if and only if $
g(u-1 )\equiv 0\pmod 2$, $gu\equiv 0\pmod k$, $k$ is even when
$u=2$, and $(k,u,g)\not\in\{(3,3,2),(3,6,2),(3,3,6),(6,2,6)\}$.
\end{theorem}

An {\it $r$-factor} is a factor which is $r$-regular. It's obvious
that a 2-factor consists of a collection of disjoint cycles. A
{\it $2$-factorization} of a graph $H$ is a partition of $E(H)$
into 2-factors. The well-known Hamilton-Waterloo problem is the
problem of determining whether $K_v$ (for $v$ odd) or $K_v$ minus
a $1$-factor (for $v$ even) has a $2$-factorization in which there
are exactly $\alpha$ $C_m$-factors and $\beta$ $C_n$-factors. For
brevity, we generalize this problem to a general graph $H$, and
use HW$(H;m,n,\alpha,\beta)$ to denote a $2$-factorization of $H$
in which there are exactly $\alpha$ $C_m$-factors and $\beta$
$C_n$-factors. So when $H=K_v$(for $v$ odd) or $K_v$ minus a
$1$-factor (for $v$ even), an HW$(H;m,n,\alpha,\beta)$ is a
solution to the original Hamilton-Waterloo problem, denoted by
HW$(v;m,n,\alpha,\beta)$. For convenience, we denote by
HWP$(v;m,n)$ the set of $(\alpha,\beta)$ for which a solution
HW$(v;m,n,\alpha,\beta)$ exists.

It is easy to see that the necessary conditions for the existence
of an HW$(v;m,n,\alpha,\beta)$ are $m | v$ when $ \alpha >0$, $n |
v$ when $ \beta >0$ and $\alpha+\beta=\lfloor \frac{v-1}{2}
\rfloor$. When $\alpha\beta=0$, the existence of an
HW$(v;m,n,\alpha,\beta)$ has been solved completely, see
Theorem~\ref{Kv}. From now on, we suppose that
$\alpha\beta\not=0$. A lot of work has been done for small values
of $m$ and $n$, especially for $m=3$. Adams et al. \cite{ABBE}
solved the case $(m,n)=(3,5)$ when $v$ is odd
 with an exception and some possible exceptions. Danziger \cite{DQS}
and Odaba\c{s}{\scriptsize I} et al. \cite{OO} solved the case
$(m,n)=(3,4)$ with three possible exceptions. Lei et al. \cite{LF}
solved the case $(m,n)=(3,7)$ when $v$ is odd with three possible
exceptions. Asplund et al. \cite{AKK} focused on $(m,n)=(3,3x)$,
and many infinite classes of HW$(v;m,n,\alpha,\beta)$s were
constructed. There are also some known results on
HW$(v;3,v,\alpha,\beta)$, see \cite{DL,DL2,HNR,LS}. For more
results on Hamilton-Waterloo problem, the reader is refer to
\cite{BuD,BD,BDD,BDT,DPR,FH,KO,LFS,MT}.

\begin{theorem}{\rm (\cite{DQS,OO})}
\label{L34} $(\alpha,\beta)\in{\rm HWP}(v;3,4)$ if and only if
$v\equiv 0\pmod{12}$, $\alpha+\beta = \lfloor \frac{v-1}{2}
\rfloor$, except possibly for
$(v,\alpha,\beta)=(24,5,6),(24,9,2),(48,17,6)$.
\end{theorem}

\begin{theorem}{\rm (\cite{ABBE})}
\label{L35} Suppose $v\equiv 15\pmod{30}$ and
$\alpha+\beta=\frac{v-1}{2}$. Then $ (\alpha,\beta)\in{\rm
HWP}(v;3,5)$ except for $(v,\alpha,\beta)=(15,6,1)$, and except
possibly for $(\alpha,\beta)=(\frac{v-3}{2},1)$ when $v>15$.
\end{theorem}

\begin{theorem}{\rm (\cite{LF})}
\label{L37} Suppose $v\equiv 21\pmod{42}$ and
$\alpha+\beta=\frac{v-1}{2}$. Then $ (\alpha, \beta)\in$ {\rm
HWP}$(v;3,7)$, except possibly for
$(v,\alpha,\beta)=(21,2,8),(21,4,6),(21,6,4)$.
\end{theorem}

Combining the known results in Theorems~\ref{L34}-\ref{L37}, we
will prove the following main result.

\begin{theorem}
\label{TH457} For $n\in\{4,5,7\}$, $(\alpha, \beta)\in{\rm
HWP}(v;3,n)$ if and only if  $v\equiv 0\pmod{3n}$, $\alpha+
\beta=\lfloor\frac{v-1}{2}\rfloor$ and
$(v,\alpha,\beta)\not=(15,6,1)$.
\end{theorem}

\section{Constructions}



Let $\Gamma$ be a finite additive group and let $S$ be a subset of
$\Gamma \backslash \{0\}$ such that the opposite of every element
of $S$ also belongs to $S$. The {\it Cayley graph} over $\Gamma$
with connection set $S$, denoted by $Cay(\Gamma, S)$, is the graph
with vertex set $\Gamma$ and edge set $E(Cay(\Gamma, S))= \{(a,b)|
a,b \in \Gamma, a-b \in S\}$. It is quite obvious that
$Cay(\Gamma, S)=Cay(\Gamma, \pm S)$.


\begin{lemma}
\label{2p-533} Let $n \geq 3$.  If $a\in Z_n$, the order of $a$
is greater than $3$ and $(i,m)=1$, then there is a $C_m$-factorization of
$Cay(Z_n \times Z_m, \pm\{0,a,2a\}\times\{ \pm i\})$.
\end{lemma}

\noindent {\it Proof:}  Since the order of $a$ is greater than
$3$, we have $|\{0,a,-a,2a,-2a\}|=5$.  Let
$C_j=((a_{j0},b_{j0})=(0,0),(a_{j1},b_{j1}),\cdots,(a_{j,m-1},b_{j,m-1}))$,
$1\leq j\leq5$, where
\[\begin{array}{lllll}
a_{11}=a, & a_{21}=0, & a_{31}=2a, & a_{41}=-a, & a_{51}=-2a,\\
a_{12}=2a, & a_{22}=-2a, & a_{32}=a, & a_{42}=-a, & a_{52}=0,
\end{array}
\]
\hspace{2.7cm} $a_{jt}=a_{j,(t-2)}$, $t\geq 3$, \\ \vspace{5pt}
\hspace*{2.6cm} $b_{jt}=ti\pmod m$, $1\le t\le m-1$.

Since $(i,m)=1$, we know that $b_{jt}$, $0\le t\le m-1$, are all
different modulo $m$. Then each $C_j$ will generate a $C_m$-factor
by $(+1\pmod{n},-)$. Thus we  can obtain the required 5
$C_m$-factors which form  a $C_m$-factorization of $Cay(Z_n \times
Z_m, \pm\{0,a,2a\}\times\{ \pm i\})$.\qed

\begin{lemma}
\label{2p-355} Let $n \geq 3$.  If $a\in Z_n$, the order of $a$ is
greater than $2$ and $(i,m)=1$,  then there is a
$C_m$-factorization of $Cay(Z_n \times Z_m, \pm\{0,a\}\times\{ \pm
i\})$.
\end{lemma}

\noindent {\it Proof:}  Because the order of $a$ is greater than
$2$, we have $|\{0,a,-a\}|=3$.  Let
$C_j=((a_{j0},b_{j0})=(0,0),(a_{j1},b_{j1}),
\cdots,(a_{j,m-1},b_{j,m-1}))$, $1\leq j\leq3$, where
\[\begin{array}{l}
a_{11}=a,  a_{21}=0,  a_{31}=-a,\\
a_{12}=a,  a_{22}=-a,  a_{32}=0,\\
a_{jt}=a_{j,(t-2)}, t\geq 3, \\
b_{jt}=ti\pmod m, 1\le t\le m-1.
\end{array}
\]
Since $(i,m)=1$, we know that $b_{jt}$, $0\le t\le m-1$, are all
different modulo $m$. Then each $C_j$ will generate a $C_m$-factor
by $(+1\pmod{n},-)$. Thus we  can obtain the required 3
$C_m$-factors which form  a $C_m$-factorization of $Cay(Z_n \times
Z_m, \pm\{0,a\}\times\{ \pm i\})$. \qed

For our recursive constructions, we need the definition of an
incomplete Hamilton-Waterloo problem design. Suppose $G$ is a
subgraph of  a graph $H$. A {\it holey 2-factor} of $H-G$ is a
2-regular subgraph of $H$ covering all vertices except those
belonging to $G$. We will also frequently speak of a holey
$C_k$-factor to mean a holey $2$-factor whose cycles all have
length $k$. Let $v-h\equiv 0\pmod 2$. An {\it incomplete
Hamilton-Waterloo problem design} on $v$ vertices with a hole of
size $h$, denoted by IHW$(v,h;m,n,\alpha,\beta,\alpha',\beta')$,
is a cycle decomposition of $K_v-E(K_h)$ if $v$ is odd, or
$K_v-E(K_h)$ minus a 1-factor $I$ if $v$ is even, such that (1)
$\alpha+\beta=\frac{v-h}{2}$,
$\alpha'+\beta'=\lfloor\frac{h-1}{2}\rfloor$; (2) there are
$\alpha$ $C_m$-factors and $\beta$ $C_n$-factors of $K_v$; (3)
there are $\alpha'$ holey $C_m$-factors and $\beta'$ holey
$C_n$-factors of $K_v-K_h$. We denote by IHWP$(v,h;m,n)$ the set
of $(\alpha,\beta,\alpha',\beta')$ for which an
IHW$(v,h;m,n,\alpha,\beta,\alpha',\beta')$ exists.

\begin{lemma}
\label{L35-IHW} $(15,0,6,1)\in {\rm IHWP}(45,15;3,5)$.
\end{lemma}

\noindent {\it Proof:}   Let the vertex set be $(Z_{6}\times Z_5)
\cup \{\infty _{i=0} ^{14}\}$.  A holey $C_5$-factor is
$Cay(Z_{6}\times Z_5, \{0\}\times \{\pm 2\})$. The required $15$
$C_3$-factors will be generated from three initial $C_3$-factors
$P_i$ $(i=1,2,3)$ by $(-, +1\pmod{5})$. For the required  six
holey $C_3$-factors, five of which can be generated from an
initial holey $C_3$-factor $Q$ by $(-, +1\pmod{5})$. The last
holey $C_3$-factor can be generated from two base cycles $(0_0,
3_4,2_0)$ and $(1_1, 5_1, 4_2)$ by $(-, +1\pmod{5})$. The cycles
of $P_i$ and $Q$ are listed below.

{\footnotesize \vspace{5pt}\noindent
\begin{tabular}{lllllllllll}
$P_1$ &$( \infty_ 0, 4_ 4, 5_ 0)$ &$( \infty_ 1, 1_ 1, 3_ 3)$ &$(
\infty_ 2, 0_ 0, 2_ 2)$ &$( \infty_ 3, 4_ 0, 5_ 2)$ &$( \infty_ 4,
2_ 3, 3_ 0)$
&$( \infty_ 5, 5_ 1, 4_ 1)$\\
&$( \infty_ 6, 3_ 4, 2_ 4)$ &$( \infty_ 7, 0_ 2, 1_ 0)$ &$(
\infty_ 8, 1_ 2, 5_ 3)$ &$( \infty_ 9, 1_ 4, 0_ 4)$ &$( \infty_
{10}, 4_ 2, 3_ 2)$
&$( \infty_ {11}, 0_ 3, 4_ 3)$\\
&$( \infty_ {12}, 0_ 1, 3_ 1)$ &$( \infty_ {13}, 2_ 0, 5_ 4)$
&$( \infty_ {14}, 1_ 3, 2_ 1)$\\
$P_2$ &$( \infty_ 0, 0_ 0, 1_ 1)$ &$( \infty_ 1, 2_ 2, 4_ 4)$ &$(
\infty_ 2, 3_ 3, 5_ 0)$ &$( \infty_ 3, 0_ 1, 3_ 4)$ &$( \infty_ 4,
1_ 2, 5_ 1)$
&$( \infty_ 5, 2_ 3, 0_ 2)$\\
&$( \infty_ 6, 4_ 0, 1_ 3)$ &$( \infty_ 7, 2_ 4, 5_ 2)$ &$(
\infty_ 8, 3_ 0, 4_ 2)$ &$( \infty_ 9, 4_ 1, 5_ 4)$ &$( \infty_
{10}, 0_ 3, 1_ 0)$
&$( \infty_ {11}, 1_ 4, 2_ 1)$\\
&$( \infty_ {12}, 2_ 0, 4_ 3)$ &$( \infty_ {13}, 3_ 1, 0_ 4)$
&$( \infty_ {14}, 5_ 3, 3_ 2)$\\
$P_3$ &$( \infty_ 0, 2_ 2, 3_ 3)$ &$( \infty_ 1, 0_ 0, 5_ 0)$ &$(
\infty_ 2, 1_ 1, 4_ 4)$ &$( \infty_ 3, 1_ 2, 2_ 3)$ &$( \infty_ 4,
0_ 1, 4_ 0)$
&$( \infty_ 5, 3_ 4, 1_ 3)$\\
&$( \infty_ 6, 5_ 1, 0_ 2)$ &$( \infty_ 7, 3_ 0, 4_ 1)$ &$(
\infty_ 8, 0_ 3, 2_ 1)$ &$( \infty_ 9, 2_ 4, 3_ 2)$ &$( \infty_
{10}, 5_ 2, 2_ 0)$
&$( \infty_ {11}, 3_ 1, 5_ 4)$\\
&$( \infty_ {12}, 5_ 3, 1_ 0)$ &$( \infty_ {13}, 1_ 4, 4_ 3)$
&$( \infty_ {14}, 4_ 2, 0_ 4)$\\
$Q$ & $( 0_ 0, 0_ 1, 4_ 2)$ &$( 1_ 1, 4_ 1, 2_ 0)$ &$( 2_ 2, 5_
2, 2_ 1)$ &$( 3_ 3, 0_ 2, 5_ 3)$ &$( 4_ 4, 2_ 4, 4_ 3)$
&$( 5_ 0, 5_ 1, 0_ 3)$\\
&$( 1_ 2, 3_ 0, 5_ 4)$ &$( 2_ 3, 1_ 3, 0_ 4)$ &$( 3_ 4, 1_ 4, 1_
0)$ &$( 4_ 0, 3_ 1, 3_ 2)$
\end{tabular}}

\qed

For next recursive construction, we still need the definition of a
cycle frame. Let $H$ be a graph $K_u[g]$ with $u$ parts
$G_1,G_2,\ldots,G_u$. A partition of $E(H)$ into holey $2$-factors
of $H-G_i(1\le i\le u)$ is said to be a {\it cycle frame of type
$g^u$}. Further, if all holey 2-factors of a cycle frame of type
$g^u$ are $C_k$-factors, then we denote the cycle frame by
$k$-CF$(g^u)$.

\begin{theorem}{\rm (\cite{S})}
\label{CF}
 There exists a $3$-{\rm CF}$(g^{u})$ if
and only if $ g\equiv 0\pmod 2$, $g(u-1)\equiv 0\pmod 3$ and
$u\geq 4$.
\end{theorem}

It's obvious that there are exactly $\frac{g}{2}$ holey 2-factors
with respect to each part. We use CF$(g^u;m,n,\alpha,\beta)$ with
$\alpha+\beta= \frac{g}{2} $ to denote a cycle frame of type $g^u$
in which there are exactly $\alpha$
 holey $C_m$-factors and $\beta$ holey $C_n$-factors with respect to each
 part. Now we use cycle frames and incomplete Hamilton-Waterloo problem
designs to give the ``Filling in Holes'' construction.

\begin{construction}
\label{C-IHW} Let  $\alpha+\beta=\frac{g}{2}$,
$\alpha'+\beta'=\lfloor\frac{h-1}{2}\rfloor$. If there exist a
{\rm CF}$(g^u;m,n,\alpha,\beta)$, an {\rm IHW}$(g+h,h
;m,n,\alpha,\beta,\alpha',\beta')$ and an {\rm HW}$(g+h;m,n,
\alpha+\alpha',\beta+\beta')$, then an {\rm HW}$(gu+h;m,n, \alpha
u+\alpha',\beta u+\beta')$ exists.
\end{construction}

\noindent {\it Proof:} We start with a CF$(g^u;m,n,\alpha,\beta)$,
for each part $G_i$, $1\le i\le u$, denote its $\alpha$ holey
$C_m$-factors by $P_{ij}(1\le j\le \alpha)$, and denote its
$\beta$ holey $C_n$-factors by $Q_{ij}(1\le j\le \beta)$.

For each $i(1\le i\le u-1)$, place a copy of an
IHW$(g+h,h;m,n,\alpha,\beta,\alpha',\beta')$ on the vertices of
the part $G_i$ and $h$ new common vertices(take the subgraph on
these $h$ vertices as the hole), whose $\alpha$ $C_m$-factors and
$\beta$ $C_n$-factors are denoted by $P'_{ij}(1\le j\le \alpha)$
and $Q'_{ij}(1\le j\le\beta)$ respectively, $\alpha'$ holey
$C_m$-factors and $\beta'$ holey $C_n$-factors are denoted by
$P''_{ij}(1\le j\le \alpha')$ and $Q''_{ij}(1\le j\le \beta')$
respectively. Further, if $h\equiv 0\pmod 2$, then $g+h\equiv
0\pmod 2$(note that the existence of a CF$(g^u)$ requires $g\equiv
0\pmod 2$). Then according to the definition of an IHW, there is a
1-factor $I_i$ of the subgraph on the vertices from $G_i$.

Place on the vertices of the part $G_u$ and these $h$ common
vertices a copy of an HW$(g+h;m,n, \alpha+\alpha',\beta+\beta')$
with $\alpha+\alpha'$ $C_m$-factors $P'_{uj}(1\le j\le
\alpha+\alpha')$ and $\beta+\beta'$ $C_n$-factors $Q'_{uj}(1\le
j\le \beta+\beta')$. If $h\equiv 0\pmod 2$, there is a 1-factor
$I_u$.

Let $$S_{ij}={P_{ij}\cup P'_{ij}}, 1\le i\le u, 1\le j\le
\alpha,$$
 $$F_{ij}={Q_{ij}\cup Q'_{ij}}, 1\le i\le u, 1\le j\le \beta,$$
$$ S_{u,j+\alpha}={(\cup_{i=1}^{u-1}P''_{ij})\cup
P'_{u,j+\alpha}},  1 \leq j \leq \alpha', $$
$$ F_{u,j+\beta}={(\cup_{i=1}^{u-1}Q''_{ij})\cup
Q'_{u,j+\beta}}, 1\leq j\leq \beta'.  $$

 Then both $S_{ij}$ and $S_{u,j+\alpha}$ are $C_m$-factors, $F_{ij}$ and $F_{u,j+\beta}$ are
$C_n$-factors on the whole vertex set, and they form an
HW$(gu+h;m,n, \alpha u+\alpha',\beta u+\beta')$. Note that if $h
\equiv 0\pmod 2$, $I={\cup_{i=1}^{u}I_i}$ is a 1-factor on the
whole vertex set. \qed

For the next recursive construction, we need more notations. When
$g(u-1)\equiv 1\pmod 2$, by Theorem~\ref{Kv} it is easy to see
that an HW$( K_u[g]; m,n, \alpha,  \beta)$ can not exist.
 In this case, by simple computation, we know that it is possible to
 partition $E(K_u[g])$ into a 1-factor, $\alpha$ $C_m$-factors and $\beta$  $C_n$-factors, where $\alpha+\beta=\lfloor \frac{g(u-1)}{2} \rfloor$.
For brevity, we still use HW$( K_u[g]; m,n, \alpha,  \beta)$  to denote such a decomposition.

\begin{construction}
\label{C-RGDD} Suppose there exist an {\rm
HW}$( K_u[g] ; m ,n, \alpha,  \beta)$ and an {\rm
HW}$(g;m,n, \alpha', \beta' )$, then an {\rm
HW}$(gu;m,n,  \alpha + \alpha',  \beta +\beta' )$ exists.
\end{construction}

\noindent {\it Proof:} We start with an HW$( K_u[g]; m,n, \alpha,
\beta)$ whose $\alpha$ $C_m$-factors are denoted by $P_{j}(1\le
j\le \alpha)$,  $\beta$ $C_n$-factors are denoted by $Q_{j}(1\le
j\le\beta)$, and a 1-factor(when $g(u-1)\equiv 1\pmod 2$) is
denoted by $I$.

For each $i(1\le i\le u)$, place a copy of an HW$(g;m,n, \alpha',
\beta' )$ on the vertices of the part $G_i$ whose $\alpha'$
$C_m$-factors and $\beta'$ $C_n$-factors are denoted by
$P'_{ij}(1\le j\le \alpha')$ and $Q'_{ij}(1\le j\le\beta')$
respectively, and a 1-factor is denoted by $I_i$ if $g\equiv
0\pmod 2$. Let $S_{j}=\bigcup_{i=1}^{u}P'_{ij}$ $(1\le j\le
\alpha')$ and $F_{j}=\bigcup_{i=1}^{u}Q'_{ij}$ $(1\le
j\le\beta')$. Then $S_{j}$ is a $C_m$-factor and $F_{j}$ is a
$C_n$-factor of the required HW$(gu;m,n,  \alpha + \alpha', \beta
+\beta')$. So we have obtained $\alpha+\alpha'$ $C_m$-factors and
$\beta+\beta'$ $C_n$-factors. At last, ${\bigcup_{i=1}^{u}I_i}$ is
a 1-factor if $g \equiv 0\pmod 2$ and $I$ is a 1-factor if $g
\equiv 1\pmod 2$ and $u \equiv 0\pmod 2$. \qed

For the next construction, we need the definition of lexicographic
product of two graphs.  Given a graph $G$, $G[n]$ is the {\it
lexicographic product} of $G$ with the empty graph on $n$ points.
Specifically, the vertex set is $\{x_i:x\in V(G), i\in Z_n\}$ and
$x_iy_j\in E(G[n])$ if and only if $xy \in E(G),i,j \in Z_n$. In
the following we will denote by $C_m[n]$ the lexicographic product
of $C_m$ with the empty graph on $n$ points.

\begin{construction}
\label{L351} If $(\alpha,\beta)\in {\rm HWP}(K_u[g];m,n)$,
$(t_i,s-t_i)\in {\rm HWP}(C_m[s];m',n')$, $1\le i\le \alpha$, and
$(r_j,s-r_j)\in {\rm HWP}(C_n[s];m',n')$, $1\le j\le \beta$, then
$(\alpha',\beta')\in {\rm HWP}(K_u[gs];m',n')$, where
$\alpha'=\sum^{\alpha}_{i=1} t_i +\sum^{\beta}_{j=1} r_j$ and
$\beta'=(\alpha+\beta)s-\alpha'$.
\end{construction}

\noindent {\it Proof:} We start with an
HW$(K_u[g];m,n,\alpha,\beta)$ with $\alpha$ $C_m$-factors and
$\beta$ $C_n$-factors. Give each vertex weight $s$, then we obtain
$\alpha$ $C_m[s]$-factors and $\beta$ $C_n[s]$-factors. Now we
replace each $C_m[s]$ and each $C_n[s]$ in the $i$-th
$C_m[s]$-factor and the $j$-th $C_n[s]$-factor with an
HW$(C_m[s];m',n',t_i,s-t_i)$ and an HW$(C_n[s];m',n',r_j,s-r_j)$
respectively. Further, take one of the $t_i$ $C_{m'}$-factors from
each HW$(C_m[s];m',n',t_i,s-t_i)$ in the $i$-th $C_m[s]$-factor,
and put them together to get a $C_{m'}$-factor of $K_u[gs]$. Thus,
we have obtained $\sum^{\alpha}_{i=1} t_i$ $C_{m'}$-factor of
$K_u[gs]$. Similarly, we can get $\sum^{\beta}_{j=1} r_j$
$C_{m'}$-factor of $K_u[gs]$ from the known
HW$(C_n[s];m',n',r_j,s-r_j)$, and $\sum^{\alpha}_{i=1}
(s-t_i)+\sum^{\beta}_{j=1}(s-r_j)=(\alpha+\beta)s-\alpha'$
$C_{n'}$-factors of $K_u[gs]$.\qed

\section{HWP$(v;3,4)$}

In this section, we will give three direct constructions and
complete the spectrum for an HW$(v;3,4,\alpha,\beta)$.

\begin{lemma}
\label{E24-34-92}   $(9,2)\in {\rm HWP}(24;3,4)$.
\end{lemma}

\noindent {\it Proof:} Let the vertex set be $\Gamma= Z_8\times
Z_3$, and the $1$-factor be $Cay(\Gamma, \{4\}\times \{0\})$.
 For the 9 $C_3$-factors,
let $F=\{Q,Q+4_0\}$, where $Q=\{( 0_0, 6_0, 0_1), ( 1_1, 1_0,
3_1),( 2_2, 2_1, 7_0),( 5_2, 3_2, 4_2)\}$.
 It's easy to see that $F$ is a $C_3$-factor
since all these 4 elements having the same subscript in $Q$ are
different modulo 4. Then $F,F+4_1,F+0_2$ are 3 $C_3$-factors. The
other $6$ $C_3$-factors can be generated from an initial
$C_3$-factor $P=\{( 0_0, 1_1, 2_2),( 3_0, 4_1, 7_1), ( 5_2, 0_2,
4_0),( 6_0, 1_0, 5_1),( 2_1, 3_2, 3_1),$ $( 6_2, 2_0, 4_2),( 7_0,
1_2, 6_1),( 0_1, 5_0, 7_2)\}$ by $(+4\pmod{8}, +1\pmod{3})$.

The required two $C_4$-factors can be generated from two base
4-cycles $( 0_0, 2_1, 5_0, 3_2)$ and $( 0_0, 5_1, 6_1, 3_1)$ by
$(+4\pmod{8}, +1\pmod{3})$ since the first coordinate of the four
elements in each cycle are different modulo 4. \qed

\begin{lemma}
\label{E24-34-56} $(5,6)\in {\rm HWP}(24;3,4)$.
\end{lemma}

\noindent {\it Proof:} Let $\Gamma=Z_8\times Z_3$. Firstly, we
construct an HW$( K_3[8]; 3, 4, 5, 3)$ with three parts $Z_8\times
\{i\}$, $i\in Z_3$. The required $5$ $C_3$-factors come from a
$C_3$-factorization of $Cay(\Gamma, \pm \{0,1,2\}\times \{\pm
1\})$ by Lemma~\ref{2p-533}. The required  three $C_4$-factors
will be generated from an initial $C_4$-factor $P=\{( 0_0, 4_1,
0_2, 5_1), ( 1_1, 4_0, 7_1, 4_2),$ $( 2_2, 6_0, 1_2, 5_0), (3_0,
0_1, 3_2, 6_1), ( 5_2, 1_0, 6_2, 2_0), ( 2_1, 7_0, 3_1, 7_2)\}$ by
$(-, +1\pmod 3)$. Then we use Construction~\ref{C-RGDD} with an
HW$(8;3,4, 0, 3)$ from Theorem~\ref{Kv} and an HW$( K_3[8]; 3, 4,
5, 3)$ constructed above to get an HW$(24; 3, 4, 5, 6)$. \qed

\begin{lemma}
\label{E48-34-176}   $(17,6)\in {\rm HWP}(48;3,4)$.
\end{lemma}
\noindent {\it Proof:}   Let the vertex set be
$\Gamma=Z_{16}\times Z_3$, and the $1$-factor be $Cay(\Gamma,
\{8\}\times \{0\})$. The required $5$ of 17 $C_3$-factors come
from a $C_3$-factorization of  $Cay(\Gamma, \pm \{0,1,2\}\times
\{\pm 1\})$ by Lemma~\ref{2p-533}. The other $12$ $C_3$-factors
can be generated from an initial $C_3$-factor $P$ by
$(+4\pmod{16}, +1\pmod{3})$. For the required 6 $C_4$-factors,
start with a cycle set $Q$ in which all these 8 elements having
the same subscript are different modulo 8. Let $F=\{Q,Q+8_0\}$.
Then $F$, $F+4_1$, $F+8_2$, $F+12_0$, $F+0_1$, $F+4_2$ are 6
$C_4$-factors. The cycles of $P$ and $Q$ are listed below.

{\footnotesize \vspace{5pt}\noindent
\begin{tabular}{lllllllllll}
$P$ & $( 0_ 0, 3_ 0, 6_ 0)$ &$( 1_ 1, 4_ 1, 8_ 2)$ &$( 2_ 2, 5_ 2,
9_ 0)$ &$( 7_ 1,11_ 2, 0_ 1)$ &$(10_ 1,14_ 2, 3_ 1)$
&$(12_ 0, 1_ 2, 6_ 1)$\\
&$(13_ 1, 5_ 0,15_ 1)$
&$(15_ 0, 7_ 2, 1_ 0)$
&$( 2_ 0,14_ 0, 8_ 1)$
&$( 4_ 2, 0_ 2, 9_ 2)$
&$( 8_ 0, 5_ 1,11_ 1)$
&$( 9_ 1, 4_ 0,15_ 2)$\\
&$(10_ 2, 7_ 0,13_ 0)$
&$(11_ 0, 2_ 1,10_ 0)$
&$(12_ 1, 3_ 2,14_ 1)$
&$(13_ 2, 6_ 2,12_ 2)$\\
$Q$ & $( 0_ 0, 8_ 2, 3_ 2,15_ 0)$ &$( 1_ 1,14_ 2, 8_ 1, 4_ 2)$
&$( 2_ 2, 1_ 2,11_ 1,15_ 1)$\\
&$( 3_ 0,12_ 1, 5_ 2,14_ 0)$
&$( 9_ 0, 5_ 0,10_ 1, 5_ 1)$
&$(12_ 0,10_ 0, 6_ 1,15_ 2)$
\end{tabular}}

\qed

Combining Theorem~\ref{L34}, Lemmas~\ref{E24-34-92}, \ref{E24-34-56} and \ref{E48-34-176},
we have the following theorem.

\begin{theorem}
\label{TH34} $(\alpha, \beta)\in {\rm HWP}(v;3,4)$ if and only if
$v\equiv 0\pmod{12}$ and $\alpha+
\beta=\lfloor\frac{v-1}{2}\rfloor$.
\end{theorem}

\section{HWP$(v;3,5)$}

In this section, we shall solve the left infinite class in
\cite{ABBE} for the existence of an HW$(v;3,5,\alpha,\beta)$ when
$v\equiv 15\pmod{30}$. Then we continue to consider the
existence of an HW$(v;3,5,\alpha,\beta)$ when $v\equiv
0\pmod{30}$.

\begin{lemma}
\label{L35-45} $(21,1)\in {\rm HWP}(45;3,5)$.
\end{lemma}

\noindent {\it Proof:}   Let the vertex set be $\Gamma
=Z_{9}\times Z_5$. The required $C_5$-factor is $Cay(\Gamma,
\{0\}\times \{\pm 1\})$. For the required $21$ $C_3$-factors, $15$
of which will be generated from an initial $C_3$-factor $P$ by
$(+3 \pmod{9}, +1 \pmod{5})$. Each cycle in $Q$ will generate a
$C_3$-factor by $(+3 \pmod{9}, +1 \pmod{5})$ since the 3 elements
in the first coordinate are different modulo 3.  Thus we have
obtained the last $6$ $C_3$-factors.  The cycles of $P$ and $Q$
are listed below.

{\footnotesize \vspace{5pt}\noindent
\begin{tabular}{lllllllllll}
$P$ &$(  1_1,  4_1,  7_4)$ &$(  4_2,  0_1,  3_4)$ &$(  7_2,  8_2,
7_0)$ &$(  8_3,  0_2,  4_0)$ &$(  6_1,  6_3,  7_3)$ &$(  2_2,
5_2,  8_4)$ &$(  3_2,  7_1,  5_1)$
&$(  3_3,  1_0,  5_4)$\\
&$(  0_4,  3_1,  8_1)$
&$(  6_4,  3_0,  2_3)$
&$(  4_4,  5_3,  1_3)$
&$(  4_3,  1_4,  2_0)$
&$(  0_0,  6_0,  8_0)$
&$(  5_0,  2_1,  2_4)$
&$(  0_3,  1_2,  6_2)$\\
$Q$ &$(  0_0,  1_1,  5_2)$ &$(  0_0,  2_2,  7_0)$ &$(  0_0,  7_1,
8_4)$ &$(  0_0,  4_2,  2_3)$ &$(  0_0,  5_3,  7_4)$ &$(  0_0,
5_1,  7_3)$
\end{tabular}}

 \qed

\begin{lemma}
\label{L35-75-105} $(\frac{v-3}{2},1)\in {\rm HWP}(v;3,5)$ for
$v=75,105$.
\end{lemma}

\noindent {\it Proof:}   Let $v=3u$ and the vertex set be
$\Gamma=Z_{u}\times Z_3$. The required $C_5$-factor is
$Cay(\Gamma, \{\pm 10\}\times \{0\})$ when $v=75$ or $Cay(\Gamma,
\{\pm 7\}\times \{0\})$ when $v=105$.

For the required $\frac{3(u-1)}{2}$ $C_3$-factors, $u$ of which
will be generated from an initial $C_3$-factor $P$ by $(+1
\pmod{u}, -)$. The other $\frac{u-3}{2}$ $C_3$-factors will be
obtained form $\frac{u-3}{2}$ 3-cycles in $Q$. Each cycle of $Q$
will generate a $C_3$-factor by $(+1 \pmod{u}, -)$
 since the first coordinate of those 3 elements of the cycle are different modulo 3.
 The cycles of $P$ and $Q$ for each $v$ are listed below.

\noindent {$v=75:$}

{\footnotesize \vspace{5pt}\noindent
\begin{tabular}{lllllllllll}
$P$
&$(  9_0, 18_0,  5_0)$
&$(  6_1, 13_2,  0_2)$
&$(  4_1, 17_0, 18_1)$
&$( 22_1,  5_1,  7_0)$
&$(  5_2, 23_0,  3_2)$
&$(  0_0, 10_1, 19_0)$\\
&$( 17_2, 12_2, 18_2)$
&$(  1_1,  8_0, 13_0)$
&$(  2_2, 19_2, 24_1)$
&$(  6_0,  8_2,  4_0)$
&$( 13_1, 19_1,  4_2)$
&$(  0_1,  9_1, 21_1)$\\
&$(  2_0,  2_1, 24_2)$
&$(  7_1, 20_2,  1_0)$
&$(  3_1, 11_0, 21_2)$
&$( 11_2, 14_2, 16_0)$
&$(  1_2, 22_2, 15_2)$
&$( 21_0, 14_0, 22_0)$\\
&$( 23_2, 12_1, 11_1)$
&$( 16_1,  6_2, 23_1)$
&$( 12_0, 15_0, 14_1)$
&$(  7_2, 16_2, 20_0)$
&$(  3_0,  8_1,  9_2)$
&$( 24_0, 10_2, 10_0)$\\
&$( 15_1, 17_1, 20_1)$\\
$Q$
&$(  0_0,  4_1,  8_2)$
&$(  0_0,  7_1, 16_2)$
&$(  0_0, 14_2,  8_1)$
&$(  0_0, 17_2,  3_1)$
&$(  0_0, 19_1, 15_2)$
&$(  0_0, 22_1, 24_2)$\\
&$(  0_0,  1_2, 21_1)$
&$(  0_0,  9_1,  9_2)$
&$(  0_0, 13_2, 14_1)$
&$(  0_0,  3_2, 11_1)$
&$(  0_0, 18_2, 20_1)$\\
\end{tabular}}

\noindent {$v=105:$}

{\footnotesize \vspace{5pt}\noindent
\begin{tabular}{lllllllllll}
$P$
& $( 3_ 1, 29_ 1, 28_ 2)$
&$( 9_ 0, 10_0,  26_0)$
&$( 1_ 0, 24_ 0, 28_ 1)$
&$( 5_ 1, 24_ 1, 26_ 1)$
&$( 0_ 2, 2_ 2, 29_ 2)$
&$( 21_ 0, 31_ 0, 9_ 1)$\\

& $( 30_ 0, 32_ 0, 32_ 1)$
&$( 2_ 1, 19_1,  17_2)$
&$( 27_ 0, 25_ 2, 34_ 2)$
&$( 10_ 1, 30_ 1, 14_ 2)$
&$(0_ 1, 6_ 1, 16_ 2)$
&$( 6_ 0, 11_ 2, 12_ 2)$\\

& $( 7_ 1, 18_ 1, 3_ 2)$
&$( 4_ 1, 34_1,  1_2)$
&$( 15_ 0, 29_ 0, 14_ 1)$
&$( 4_ 2, 23_ 2, 27_ 2)$
&$( 0_ 0, 10_ 2, 15_ 2)$
&$( 8_ 0, 19_ 0, 23_ 0)$\\

& $( 14_ 0, 5_ 2, 22_ 2)$
&$( 3_ 0, 25_0,  15_1)$
&$( 4_ 0, 7_ 0, 13_ 0)$
&$( 18_ 2, 21_ 2, 32_ 2)$
&$( 2_ 0, 20_ 2, 30_ 2)$
&$( 1_ 1, 23_ 1, 33_ 1)$\\

& $( 12_ 0, 20_ 0, 13_ 1)$
&$( 8_ 1, 12_ 1, 20_ 1)$
&$( 28_ 0, 33_0,  17_1)$
&$( 21_ 1, 22_ 1, 8_ 2)$
&$( 16_ 0, 13_ 2, 33_ 2)$
&$( 22_ 0, 9_ 2, 31_ 2)$\\

& $( 5_ 0, 11_ 1, 24_ 2)$
&$( 11_ 0, 16_ 1, 7_ 2)$
&$( 17_ 0, 25_ 1, 19_ 2)$
&$( 18_ 0, 27_ 1, 6_ 2)$
&$( 34_ 0, 31_ 1, 26_ 2)$\\
$Q$
&$( 0_ 0, 3_ 1, 3_ 2)$
&$( 0_ 0, 7_ 1, 0_ 2)$
&$( 0_ 0, 10_ 1, 11_ 2)$
&$( 0_ 0, 11_ 1, 16_ 2)$
&$( 0_ 0, 14_ 1, 21_ 2)$
&$( 0_ 0, 15_ 1, 4_ 2)$\\
&$( 0_ 0, 16_ 1, 24_ 2)$
&$( 0_ 0, 17_ 1, 20_ 2)$
&$( 0_ 0, 18_ 1, 29_ 2)$
&$( 0_ 0, 21_ 1, 30_ 2)$
&$( 0_ 0, 22_ 1, 34_ 2)$
&$( 0_ 0, 26_ 1, 14_ 2)$\\
&$( 0_ 0, 29_ 1, 12_ 2)$
&$( 0_ 0, 30_ 1, 1_ 2)$
&$( 0_ 0, 31_ 1, 13_ 2)$
&$( 0_ 0, 33_ 1, 25_ 2)$
\end{tabular}}

\qed

\begin{lemma}
\label{L35odd}  If $v \equiv 15 \pmod {30}$ and $v>15$, then
$(\frac{v-3}{2},1)\in {\rm HWP}(v;3,5)$.
\end{lemma}

\noindent {\it Proof:}  Let $v =30u+15$, $u>0$. For $u\le 3$, the
conclusion comes from Lemmas~\ref{L35-45} and~\ref{L35-75-105}.
Applying Construction~\ref{C-IHW} with an IHW$(45,15
;3,5,15,0,6,1)$ from Lemma~\ref{L35-IHW}, a CF$(30^u;3,5,15,0)$
from Theorem~\ref{CF} and an HW$(45;3,5,21,1)$ from
Lemma~\ref{L35-45}, we get an HW$(v;3,5, \frac{v-3}{2},1)$ for any
$u\ge 4$. \qed

\begin{lemma}
\label{L35-30}   $ (\alpha,\beta)\in {\rm HWP}(30;3,5)$ if and
only if $\alpha+\beta=14$.
\end{lemma}

\noindent {\it Proof:} Let the vertex set be $\Gamma= Z_{10}\times
Z_3$ and the $1$-factor be $Cay(\Gamma, \{5\}\times \{0\})$. For
$(\alpha,\beta)=(10,4)$, we get the conclusion  by using
Construction~\ref{C-RGDD} with an HW$(10;3,5, 0, 4 )$ and an HW$(
K_3[10] ; 3,5, 10,  0)$ from Theorem~\ref{Kv}. For all the other
cases, the methods of generating the required $\alpha$
$C_3$-factors and $\beta$ $C_5$-factors are listed in Table 1. For
$C_3$-factors, here are three methods. (1) From a
$C_3$-factorization of certain Cayley graphs; (2) From several
initial $C_3$-factors $P_i$s by $(-, +1\pmod 3)$ or $(+2\pmod
{10}, -)$; (3) From several cycle sets $Q_i$s by $(+2\pmod {10},
-)$, note that each $Q_i$ will generate a $C_3$-factor by
$(+2\pmod {10}, -)$, since the two elements having the same
subscript in $Q_i$ have different parity. For $C_5$-factors, only
the first two methods are applied, and the initial $C_5$-factors
are denoted by $P_i'$ in Table 1. For the sake of brevity, we list
the cycles of $P_i$, $P_i'$ and $Q_i$ in Appendix A.\qed

\begin{lemma}
\label{RGDD60} For each $(\alpha,\beta)\in
\{(0,22),(6,16),(12,10)\}$, $(\alpha,\beta)\in {\rm
HWP}(K_4[15];3,5)$.
\end{lemma}

\noindent {\it Proof:} Let the vertex set be $Z_{15}\times Z_4$,
and the four parts of $K_4[15]$ be $Z_{15}\times \{i\}$, $i\in
Z_4$. The $\alpha$ $C_3$-factors will be obtained from $\alpha$
3-cycles from a cycle set $T$ by $(+3 \pmod {15},+1 \pmod {4})$.
 Note that the first coordinate of the 3 elements
in each cycle from $T$ are different modulo 3, so each cycle of
$T$ will generate a $C_3$-factor by $(+3 \pmod {15}, +1 \pmod
{4})$.

{\scriptsize
\begin{center}

\begin{tabular}{|c|l|l|l|}
\multicolumn{3}{c}{\bf Table 1 HWP$(30;3,5)$}\\[5pt]
  \hline
 $(\alpha, \beta)$  & \hspace*{1.5cm} $C_3$-factor &  $\hspace*{1.5cm} C_5$-factor \\   \hline

 $(1,13)$  &  \hspace*{2pt}  1: $Cay(\Gamma, \{0\}\times \{\pm 1\})$ &  12: $P_1'$, $P_2'$, $P_3'$, $P_4'$ $(-, +1\pmod 3)$\\
     & &                                 \hspace*{2pt} 1: $Cay(\Gamma, \{\pm 2\}\times \{0\})$       \\
 \hline

                                                 & & 10: $P_1'$,$P_2'$  $(+2\pmod {10}, -)$ \\
 $(2,12)$ &  \hspace*{2pt}  2: $Q_1$,$Q_2$   &    \hspace*{2pt} 1: $Cay(\Gamma, \{\pm 2\}\times \{0\})$   \\
                                                 & &  \hspace*{2pt} 1: $Cay(\Gamma, \{\pm 4\}\times \{0\})$  \\
 \hline

 $(3,11)$ &  \hspace*{2pt} 2: $Q_1$,$Q_2$        &   10: $P_1'$,$P_2'$  $(+2\pmod {10}, -)$  \\
         &  \hspace*{2pt} 1: $Cay(\Gamma, \{0\}\times \{\pm 1\})$    &   \hspace*{2pt}  1: $Cay(\Gamma, \{\pm 2\}\times \{0\})$   \\

 \hline
 $(4,10)$ &  \hspace*{2pt}  4: $Q_1$,$Q_2$,$Q_3$,$Q_4$    & 10: $P_1'$,$P_2'$  $(+2\pmod {10}, -)$    \\
 \hline

 $(5,9)$   &  \hspace*{2pt}  5:  Lemma~\ref{2p-533} with $a=i=1$                & \hspace*{2pt}  9: $P_1'$,$P_2'$,$P_3'$ $(-, +1\pmod 3)$   \\
 \hline

                                              & &  \hspace*{2pt}  6: $P_1'$,$P_2'$  $(-, +1\pmod 3)$ \\
 $(6,8)$ &  \hspace*{2pt} 6: $P_1$,$P_2$ $(-, +1\pmod 3)$    &    \hspace*{2pt}  1: $Cay(\Gamma, \{\pm 2\}\times \{0\})$   \\
                                              & &  \hspace*{2pt}  1: $Cay(\Gamma, \{\pm 4\}\times \{0\})$  \\
 \hline

 $(7,7)$ & \hspace*{2pt} 6: $P_1$,$P_2$  $(-, +1\pmod 3)$          &  \hspace*{2pt}  6: $P_1'$,$P_2'$  $(-, +1\pmod 3)$  \\
        &  \hspace*{2pt} 1: $Cay(\Gamma, \{0\}\times \{\pm 1\})$    &  \hspace*{2pt}  1: $Cay(\Gamma, \{\pm 2\}\times \{0\})$   \\

 \hline

 $(8,6)$ &  \hspace*{2pt} 5: $P_1$  $(+2\pmod {10}, -)$               &  \hspace*{2pt}  5:  $P_1'$  $(+2\pmod {10}, -)$   \\
         &  \hspace*{2pt} 3: $Q_1$,$Q_2$,$Q_3$    &  \hspace*{2pt}  1: $Cay(\Gamma, \{\pm 2\}\times \{0\})$   \\
 \hline

                                                          &  \hspace*{2pt} 5: $P_1$ $(+2\pmod {10}, -)$ & \\
 $(9,5)$ &  \hspace*{2pt}  3: $Q_1$,$Q_2$,$Q_3$      &  \hspace*{2pt} 5:  $P_1'$  $(+2\pmod {10}, -)$   \\
                                                          &  \hspace*{2pt} 1: $Cay(\Gamma, \{ 0 \}\times \{\pm 1\})$ & \\
 \hline

 $(11,3)$ &  \hspace*{2pt}  6: $P_1$,$P_2$  $(-, +1\pmod 3)$               & \hspace*{2pt}  3:  $P_1'$  $(-, +1\pmod 3)$ \\
          &  \hspace*{2pt}  5:  Lemma~\ref{2p-533} with $a=i=1$ & \\
 \hline

 $(12,2)$ & 10: $P_1$,$P_2$  $(+2\pmod {10}, -)$           &  \hspace*{2pt}  1: $Cay(\Gamma, \{\pm 2\}\times \{0\})$  \\
          &  \hspace*{2pt}  2: $Q_1$,$Q_2$             &  \hspace*{2pt}  1: $Cay(\Gamma, \{\pm 4\}\times \{0\})$   \\
 \hline

 $(13,1)$ & 10: $P_1$,$P_2$  $(+2\pmod {10}, -)$           &  \hspace*{2pt}  1: $Cay(\Gamma, \{\pm 2\}\times \{0\})$ \\
          &  \hspace*{2pt}  3: $Q_1$,$Q_2$,$Q_3$   & \\
 \hline
\end{tabular}
\end{center}
}

For $\beta$ $C_5$-factors, ten of them will be obtained from two
cycle sets $Q_1'$ and $Q_2'$. Here $\{Q_i' + (5j+k)_{0}\ |\
j=0,1,2\}$ is a $C_5$-factor for any $i=1,2$ and $k=0,1,\cdots,4$
since these 5 elements having the same subscript in $Q_i'$ are
different modulo 5. The other $\beta-10$ $C_5$-factors will be
obtained from $\beta-10$ 5-cycles in a cycle set $T'$ by $(+5
\pmod {15}, +1 \pmod {4})$, since the first coordinate of the 5
elements in each cycle from $T'$ are different modulo 5. We list
the 1-factor $I$ and the cycles in $T$, $Q_1'$, $Q_2'$ and $T'$ in
Appendix B.\qed

\begin{lemma}
\label{L35-60}  $ (\alpha,\beta)\in {\rm HWP}(60;3,5)$ if and only
if $\alpha+\beta=29$.
\end{lemma}

\noindent {\it Proof:} By Theorem~\ref{L35}, there is an
HWP$(15;3,5,\alpha_1,7-\alpha_1)$ for any $0\leq \alpha_1 \leq7$
and $\alpha_1 \neq 6$. Apply Construction~\ref{C-RGDD} with an
HW$( K_4[15];3,5,\alpha_2,22-\alpha_2)$ for $\alpha_2=0,6,12$ from
Lemma~\ref{RGDD60} to get an HW$(
60;3,5,\alpha_1+\alpha_2,29-\alpha_1-\alpha_2)$. Thus we have
$(\alpha,\beta)\in$ HWP$(60;3,5)$ for $0\le \alpha\le 19$ and
$\alpha\not=18$.

Similarly, for $(\alpha,\beta)=(20,9),(25,4)$, an
HW$(60;3,5,\alpha,\beta)$ can be obtained from the existence of an
HW$(K_3[20]; 3, 5,  20,  0)$, an HW$(20;3,5, 0, 9)$, an HW$(
K_6[10] ; 3,5, 25,  0)$ and an HW$(10;3,5, 0, 4)$ from
Theorem~\ref{Kv}.

For all the other cases, let the vertex set be $\Gamma=
Z_{15}\times Z_4$ and the $1$-factor be $Cay(\Gamma, \{0\}\times
\{2\})$, the methods of generating the required $\alpha$
 $C_3$-factors and $\beta$ $C_5$-factors are given in Table 2.
For generating $C_3$-factors, here are five methods. (1) From a
$C_3$-factorization of certain Cayley graphs; (2) From an initial
$C_3$-factor $P$ by $(+ 1\pmod {15}, -)$; (3) From several cycle
sets $Q_i$s, note that $\{Q_i + (3j+k)_{0}\ |\ j=0,1,\cdots,4\}$
is a $C_3$-factor for $k=0,1,2$ since these 3 elements having the
same subscript in $Q_i$ are different modulo 3; (4) From a cycle
in $T$ by $(+1 \pmod {15}, +1 \pmod {4})$. A $C_3$-factor $F$ can
be obtained from the cycle in $T$ by $(+3 \pmod {15}, +1 \pmod
{4})$. Then three $C_3$-factors can be generated from $F$ by $(+i
\pmod {15}, -)$, $i=0,1,2$; (5) From a cycle set $S$ by $(+3 \pmod
{15},+1 \pmod {4})$. Note that the first coordinate of the 3
elements in each cycle from $S$ are different modulo 3, so each
cycle of $S$ will generate a $C_3$-factor by $(+3 \pmod {15}, +1
\pmod {4})$.

 For $C_5$-factors, three methods are applied.
(1) From a $C_5$-factorization of certain Cayley graphs; (2) From
a cycle set $Q'$, where $\{Q' + (5j+k)_{0}\ |\ j=0,1,2\}$ is a
$C_5$-factor for $k=0,1,\cdots,4$ since these 5 elements having
the same subscript in $Q'$ are different modulo 5; (3) From a
cycle set $T'$ by $(+5 \pmod {15}, +1 \pmod {4})$. The cycles of
$P$, $Q_i$, $Q'$, $S$, $T$ and $T'$ are given in
 Appendix C. \qed

{\scriptsize
\begin{center}

\begin{tabular}{|c|l|l|l|}
\multicolumn{3}{c}{\bf Table 2 HWP$(60;3,5)$}\\[5pt]
  \hline

 $(\alpha, \beta)$  & \hspace*{1.5cm} $C_3$-factor &  $\hspace*{1.5cm} C_5$-factor \\   \hline

            &  15: $P$                & \hspace*{2pt} 5: $Q'$    \\
 $(18,11)$  & \hspace*{2pt} 3: $Q_1$                  & \hspace*{2pt} 5: $T'$   \\
                                                          &  & \hspace*{2pt} 1: $Cay(\Gamma, \{\pm 6\}\times \{0\})$   \\
 \hline

          &  15: $P$                 & \hspace*{2pt} 5: $Q'$     \\
  $(21,8)$ & \hspace*{2pt} 3: $Q_1$              & \hspace*{2pt} 3: $T'$   \\
           & \hspace*{2pt} 3: $T$         &  \\
 \hline

            & 15: $P$                  & \hspace*{2pt} 5: $Q'$    \\
  $(22,7)$  & \hspace*{2pt} 3: $Q_1$     &\hspace*{2pt}  1: $Cay(\Gamma, \{\pm 3\}\times \{0\})$   \\
            & \hspace*{2pt} 3: $T$     &\hspace*{2pt}  1: $Cay(\Gamma, \{\pm 6\}\times \{0\})$ \\
            & \hspace*{2pt} 1: $Cay(\Gamma, \{ \pm5 \} \times \{0\})$      &  \\
 \hline

             & 15: $P$           & \hspace*{2pt} 5: $Q'$    \\
  $(23,6)$ &  \hspace*{2pt} 3: $Q_1$          & \hspace*{2pt} 1: $Cay(\Gamma, \{\pm 6\}\times \{0\})$   \\
            & \hspace*{2pt} 5:  $S$      &  \\
 \hline

            &  15: $P$               & \hspace*{2pt} 3: $T'$   \\
 $(24,5)$   & \hspace*{2pt} 9: $Q_i$,$1 \leq i \leq3$     & \hspace*{2pt} 1: $Cay(\Gamma, \{\pm 3\}\times \{0\})$   \\
                                                          &  & \hspace*{2pt} 1: $Cay(\Gamma, \{\pm 6\}\times \{0\})$   \\
 \hline

          &  15: $P$               &   \\
  $(26,3)$ & \hspace*{2pt} 6: $Q_1$,$Q_2$             & \hspace*{2pt} 3: $T'$  \\
           & \hspace*{2pt} 5:  $S$     &  \\
 \hline

  $(27,2)$ & 15: $P$                    & \hspace*{2pt} 1: $Cay(\Gamma, \{\pm 3\}\times \{0\})$  \\
           & 12: $Q_i$,$1 \leq i \leq4$       & \hspace*{2pt} 1: $Cay(\Gamma, \{\pm 6\}\times \{0\})$   \\
 \hline

           &  15: $P$              &   \\
  $(28,1)$ &  12: $Q_i$,$1 \leq i \leq4$          & \hspace*{2pt} 1: $Cay(\Gamma, \{\pm 6\}\times \{0\})$   \\
           & \hspace*{2pt} 1:  $Cay(\Gamma, \{ \pm5 \} \times \{0\})$      &  \\
 \hline

\end{tabular}
\end{center}
}

\begin{lemma}
\label{L35even}  If $v \equiv 0\pmod {30}$, then
$(\alpha,\beta)\in {\rm HWP}(v;3,5)$ for any
$\alpha+\beta=\frac{v-2}{2}$.
\end{lemma}

\noindent {\it Proof:} Let $v=30u$, $u\ge 1$. For $u\le 2$, the
conclusion comes from Lemmas~\ref{L35-30} and~\ref{L35-60}. For
$u=3$, start with an HW$(K_3[3];3,5,3,0)$, an
HW$(C_3[10];3,5,10,0)$ and an HW$(C_3[10];$ $3,5,0,10)$ from
Theorem~\ref{Kv}, apply Construction~\ref{L351} with $s=10$ and
$t_i\in\{0,10\}$ to get an
HW$(K_3[30];3,5,\sum_{i=1}^{3}t_i,30-\sum_{i=1}^{3}t_i)$. Then we
apply Construction~\ref{C-RGDD} with an
HW$(30;3,5,\alpha',14-\alpha')$, $0\leq \alpha'\leq14$, from
Lemma~\ref{L35-30} to obtain an
HW$(90;3,5,\sum_{i=1}^{3}t_i+\alpha',30-\sum_{i=1}^{3}t_i+(14-\alpha'))$.
Thus we have obtained an HW$(90;3,5,\alpha,\beta)$ for any
$\alpha+\beta=44$ since $\sum_{i=1}^{3}t_i+\alpha'$ can cover the
integers from $0$ to $44$.

For $u\geq 4$, similarly, we start with an
HW$(K_u[6];3,5,3u-3,0)$, an HW$(C_3[5];3,5,5,0)$ and an
HW$(C_3[5];3,5,0,5)$ from Theorem~\ref{Kv}, and apply
Construction~\ref{L351} with $s=5$ and $t_i\in\{0,5\}$ to get an
HW$(K_u[30];3,5,\sum_{i=1}^{3u-3}t_i,15u-15-\sum_{i=1}^{3u-3}t_i)$.
Further, applying Construction~\ref{C-RGDD} with an
HW$(30;3,5,\alpha',14-\alpha')$, $0\leq \alpha'\leq14$, from
Lemma~\ref{L35-30}, we can obtain an
HW$(30u;3,5,\alpha'+\sum_{i=1}^{3u-3}t_i,14-\alpha'+15u-15-\sum_{i=1}^{3u-3}t_i)$.
It's easy to prove that $\alpha'+\sum_{i=1}^{3u-3}t_i$ can cover
the integers from $0$ to $15u-1$. The proof is complete.\qed

Combining Theorem~\ref{L35}, Lemmas~\ref{L35odd} and ~\ref{L35even},
we have the following theorem.

\begin{theorem}
\label{TH35} $ (\alpha, \beta)\in {\rm HWP}(v;3,5)$ if and only if
$v\equiv 0\pmod{15}$, $\alpha+ \beta=\lfloor\frac{v-1}{2}\rfloor$
and $(\alpha,\beta,v)\not=(6,1,15)$.
\end{theorem}

\section{HWP$(v;3,7)$}

In this section, we shall solve the three left cases in \cite{LF}
for the existence of an HW$(v;3,7,\alpha,\beta)$ when $v\equiv
21\pmod{42}$. Then we continue to consider  the case $v\equiv
0\pmod{42}$.

\begin{lemma}
\label{E21-37} For each  $(\alpha,\beta)\in\{(2,8),(4,6),(6,4)\}$,
$(\alpha,\beta)\in {\rm HWP}(21;3,7)$.
\end{lemma}

\noindent {\it Proof:} Let the vertex set be $\Gamma=Z_{7}\times
Z_{3}$.  For $\alpha =2$, the required two $C_{3}$-factors can be
 generated from two base cycles $(0_0, 1_1, 2_2)$ and
 $(0_0, 4_1, 1_2)$ by $(+1\pmod 7 ,-)$. Seven of the required $C_{7}$-factors can be obtained from
 an initial $C_{7}$-factor $P=\{(0_0, 5_2, 6_0, 1_1, 3_0,  3_1, 3_2),$ $(2_2,  5_0, 4_1, 6_2, 4_2, 5_1, 6_1),(0_1, 2_1, 0_2, 1_2, 1_0, 2_0, 4_0)\}$
 by $(+1\pmod 7 ,-)$. The last $C_{7}$-factor is $Cay(\Gamma,
  \{\pm3\}\times \{0\})$.

For $\alpha =4$, a $C_3$-factor is $Cay(\Gamma, \{0\}\times \{\pm
1\})$, and the other 3 $C_3$-factors can be generated from an
initial $C_3$-factor $P$ by $(-,  +1\pmod 3)$.  All $C_7$-factors
can be generated from two initial $C_7$-factors $Q_1$ and $Q_2$ by
$(-, +1\pmod 3)$. $P$, $Q_1$ and $Q_2$ are listed below.

{\footnotesize \noindent
\begin{tabular}{llllllllll}
$P$ & $( 0_0, 1_1, 2_2)  $ & $( 3_0, 4_1, 5_2)  $ & $( 6_0, 0_1,
3_1)  $ & $( 1_2, 4_2, 6_1)  $ & $( 2_0, 0_2, 4_0)  $ & $( 5_0,
1_0, 3_2)  $ & $( 2_1, 5_1, 6_2)  $
\end{tabular}}

{\footnotesize \noindent
\begin{tabular}{lllllllll}
$Q_1$ & $( 0_0, 5_2, 1_1, 0_1, 3_0, 2_2, 1_2)  $ & & $( 4_1, 2_0,
6_0, 5_0, 0_2, 6_2, 3_2)  $ & & $( 3_1, 2_1, 4_0, 6_1, 1_0, 4_2,
5_1)  $ \\
$Q_2$ & $( 0_0, 2_0, 1_1, 4_2, 6_2, 0_1, 4_0) $ & & $( 2_2, 5_0,
2_1, 3_0, 0_2, 5_2, 6_1)  $ & & $( 4_1, 3_1, 1_0, 6_0, 3_2, 1_2,
5_1)  $\vspace{5pt}
\end{tabular}}

For $\alpha =6$, All $C_3$-factors can be generated from two
initial $C_3$-factors $P_1$ and $P_2$ by $(-, +1\pmod 3)$. Three
$C_7$-factors can be generated from an initial $C_7$-factor $Q$ by
$(-, +1\pmod 3)$. The last $C_{7}$-factor is $Cay(\Gamma,
  \{\pm3\}\times \{0\})$. $P_1$, $P_2$ and $Q$ are listed below.

{\footnotesize \vspace{5pt}\noindent
\begin{tabular}{llllllllll}
$P_1$
& $( 0_0, 1_1, 2_2)  $
& $( 3_0, 4_1, 5_2)  $
& $( 6_0, 0_1, 4_2)  $
& $( 1_2, 3_1, 0_2)  $
& $( 2_0, 6_1, 3_2)  $
& $( 5_0, 2_1, 4_0)  $
& $( 1_0, 5_1, 6_2)  $ \\
$P_2$ & $( 0_0, 5_2, 0_1)  $ & $( 1_1, 1_2, 2_1)  $ & $( 2_2, 3_0,
3_2)  $ & $( 4_1, 3_1, 4_0)  $ & $( 6_0, 5_0, 5_1)  $ & $( 2_0,
0_2, 6_2)  $ & $( 4_2, 6_1, 1_0)  $
\end{tabular}}

{\footnotesize \noindent
\begin{tabular}{lllllllll}
$Q$
& $( 0_0, 1_2, 6_2, 6_1, 4_1, 2_1, 2_0)  $ &
& $( 1_1, 4_2, 5_1, 0_1, 4_0, 2_2, 5_0)  $ &
& $( 3_0, 0_2, 6_0, 3_2, 5_2, 3_1, 1_0)  $
\end{tabular}}

\qed

\begin{lemma}
\label{L37-42} $ (\alpha,\beta)\in {\rm HWP}(42;3,7)$ if and only
if $\alpha+\beta=20$.
\end{lemma}

\noindent {\it Proof:} For $(\alpha,\beta)=(14,6)$, we get the
conclusion  by using Construction~\ref{C-RGDD} with an HW$(14;$ $
3, 7, 0, 6)$ and an HW$( K_3[14]; 3, 7, 14,  0)$ from
Theorem~\ref{Kv}.

For all the other cases, let the vertex set be $\Gamma=
Z_{14}\times Z_3$ and the $1$-factor be $Cay(\Gamma, \{7\}\times
\{0\})$. The methods of generating the required $\alpha$
$C_3$-factors and $\beta$ $C_7$-factors are listed in the
following table.

{\scriptsize
\begin{center}

\begin{tabular}{|c|l|l|l|}
\multicolumn{3}{c}{\bf Table 3 HWP$(42;3,7)$}\\[5pt]
  \hline

  $(\alpha, \beta)$  & \hspace*{1.5cm} $C_3$-factor &  $\hspace*{2cm} C_7$-factor \\   \hline

 $(1,19)$  &\hspace*{2pt} 1: $Cay(\Gamma, \{0\}\times \{\pm 1\})$ &  18: $P_1'$,$P_2'$,$P_3'$  $(+7 \pmod {14}, +1\pmod 3)$  \\
                                                   & &  \hspace*{2pt} 1: $Cay(\Gamma, \{\pm 4\}\times \{0\})$   \\
 \hline

  $(2,18)$ &\hspace*{2pt} 2: $Q_1$,$Q_2$     &    14: $P_1'$,$P_2'$  $(+2\pmod {14}, -)$ \\
                                                   & &  \hspace*{2pt} 4: $T'$   \\

 \hline

                                                        &   &   14: $P_1'$,$P_2'$  $(+2\pmod {14}, -)$  \\
 $(3,17)$ & \hspace*{2pt} 3: $Q_i$, $1 \leq i\leq 3$          &  \hspace*{2pt} 1: $Cay(\Gamma, \{\pm 2\}\times \{0\})$   \\
                                                          &  &  \hspace*{2pt} 1: $Cay(\Gamma, \{\pm 4\}\times \{0\})$   \\
                                                          &  &  \hspace*{2pt} 1: $Cay(\Gamma, \{\pm 6\}\times \{0\})$   \\

 \hline
                                                            &  & 14: $P_1'$,$P_2'$  $(+2\pmod {14}, -)$    \\
 $(4,16)$ & \hspace*{2pt} 4: $Q_i$, $1 \leq i\leq 4$     &  \hspace*{2pt} 1: $Cay(\Gamma, \{\pm 2\}\times \{0\})$   \\
                                                           &  &  \hspace*{2pt} 1: $Cay(\Gamma, \{\pm 4\}\times \{0\})$   \\
 \hline

                                                        &   &   12: $P_1'$,$P_2'$   $(+7 \pmod {14}, +1\pmod 3)$   \\
 $(5,15)$  & \hspace*{2pt} 5:  Lemma~\ref{2p-533} with $a=2$ and $i=1$               &  \hspace*{2pt} 1: $Cay(\Gamma, \{\pm 2\}\times \{0\})$   \\
                                                          &  &  \hspace*{2pt} 1: $Cay(\Gamma, \{\pm 4\}\times \{0\})$   \\
                                                          &  &  \hspace*{2pt} 1: $Cay(\Gamma, \{\pm 6\}\times \{0\})$   \\

 \hline

                                                       & & 12: $P_1'$,$P_2'$   $(+7 \pmod {14}, +1\pmod 3)$ \\
 $(6,14)$ & \hspace*{2pt} 6: $P_1$  $(+7 \pmod {14}, +1\pmod 3)$     &   \hspace*{2pt} 1: $Cay(\Gamma, \{\pm 2\}\times \{0\})$   \\
                                                       & & \hspace*{2pt} 1: $Cay(\Gamma, \{\pm 4\}\times \{0\})$  \\
 \hline

 $(7,13)$ & \hspace*{2pt} 6: $P_1$   $(+7 \pmod {14}, +1\pmod 3)$      & 12: $P_1'$,$P_2'$   $(+7 \pmod {14}, +1\pmod 3)$  \\
          & \hspace*{2pt} 1: $Cay(\Gamma, \{0\}\times \{\pm 1\})$      & \hspace*{2pt} 1: $Cay(\Gamma, \{\pm 4\}\times \{0\})$   \\
 \hline

 $(8,12)$ & \hspace*{2pt} 5:  Lemma~\ref{2p-533} with $a=11$ and $i=1$             &  12: $P_1'$,$P_2'$   $(+7 \pmod {14}, +1\pmod 3)$ \\
          & \hspace*{2pt} 3: $T$     &   \\
 \hline

                                                                             &   &  \hspace*{2pt} 6: $P_1'$    $(+7 \pmod {14}, +1\pmod 3)$   \\
 $(9,11)$  & \hspace*{2pt} 9: $T$                   & \hspace*{2pt} 3: $Q'$     \\
                                                                               &  & \hspace*{2pt} 1: $Cay(\Gamma, \{\pm 2\}\times \{0\})$   \\
                                                                               &  & \hspace*{2pt} 1: $Cay(\Gamma, \{\pm 4\}\times \{0\})$   \\
 \hline

                                                                              &   & \hspace*{2pt} 7: $P_1'$   $(+2 \pmod {14}, -)$   \\
 $(10,10)$  & \hspace*{2pt} 5:  Lemma~\ref{2p-533} with $a=11$ and $i=1$                                   & \hspace*{2pt} 1: $Cay(\Gamma, \{\pm 2\}\times \{0\})$    \\
           &  \hspace*{2pt} 5: $Q_i$, $1 \leq i\leq 5$                     &\hspace*{2pt}  1: $Cay(\Gamma, \{\pm 4\}\times \{0\})$   \\
                                                                               &  &  \hspace*{2pt} 1: $Cay(\Gamma, \{\pm 6\}\times \{0\})$   \\
 \hline

                                                                              &   & \hspace*{2pt} 6: $P_1'$    $(+7 \pmod {14}, +1\pmod 3)$  \\
 $(11,9)$  & \hspace*{2pt} 5:  Lemma~\ref{2p-533} with $a=11$   and $i=1$                             & \hspace*{2pt} 1: $Cay(\Gamma, \{\pm 2\}\times \{0\})$    \\
           & \hspace*{2pt} 6: $P_1$   $(+7 \pmod {14}, +1\pmod 3)$                              & \hspace*{2pt} 1: $Cay(\Gamma, \{\pm 4\}\times \{0\})$   \\
                                                                               &  & \hspace*{2pt} 1: $Cay(\Gamma, \{\pm 6\}\times \{0\})$   \\
 \hline

                                                              & & \hspace*{2pt} 6: $P_1'$   $(+7 \pmod {14}, +1\pmod 3)$ \\
 $(12,8)$ & 12: $P_1$,$P_2$  $(+7 \pmod {14}, +1\pmod 3)$     &   \hspace*{2pt} 1: $Cay(\Gamma, \{\pm 2\}\times \{0\})$   \\
                                                              & & \hspace*{2pt} 1: $Cay(\Gamma, \{\pm 4\}\times \{0\})$  \\
 \hline

 $(13,7)$ & 12: $P_1$,$P_2$   $(+7 \pmod {14}, +1\pmod 3)$      & \hspace*{2pt} 6: $P_1'$   $(+7 \pmod {14}, +1\pmod 3)$  \\
          & \hspace*{2pt} 1: $Cay(\Gamma, \{0\}\times \{\pm 1\})$            & \hspace*{2pt} 1: $Cay(\Gamma, \{\pm 4\}\times \{0\})$   \\
 \hline

           & \hspace*{2pt}  6: $P_1$   $(+7 \pmod {14}, +1\pmod 3)$               & \hspace*{2pt} 3: $Q'$   \\
 $(15,5)$  & \hspace*{2pt} 9: $T$   & \hspace*{2pt} 1: $Cay(\Gamma, \{\pm 2\}\times \{0\})$   \\
                                                                 &  & \hspace*{2pt} 1: $Cay(\Gamma, \{\pm 4\}\times \{0\})$   \\
 \hline

          & \hspace*{2pt} 6: $P_1$   $(+7 \pmod {14}, +1\pmod 3)$                & \hspace*{2pt} 3: $Q'$   \\
 $(16,4)$ & \hspace*{2pt} 9: $T$     & \hspace*{2pt} 1: $Cay(\Gamma, \{\pm 4\}\times \{0\})$   \\
          & \hspace*{2pt} 1: $Cay(\Gamma, \{0\}\times \{ \pm 1\})$  & \\
 \hline

           & 14: $P_1$   $(+1 \pmod {14}, -)$              & \hspace*{2pt} 1: $Cay(\Gamma, \{\pm 2\}\times \{0\})$    \\
  $(17,3)$ & \hspace*{2pt} 2: $S$          & \hspace*{2pt} 1: $Cay(\Gamma, \{\pm 4\}\times \{0\})$   \\
           & \hspace*{2pt} 1: $Cay(\Gamma, \{0\}\times \{ \pm 1\})$      & \hspace*{2pt} 1: $Cay(\Gamma, \{\pm 6\}\times \{0\})$ \\
 \hline

            & 14: $P_1$   $(+1 \pmod {14}, -)$              & \hspace*{2pt} 1: $Cay(\Gamma, \{\pm 2\}\times \{0\})$    \\
  $(18,2)$ &  \hspace*{2pt} 3: $S$     & \hspace*{2pt} 1: $Cay(\Gamma, \{\pm 4\}\times \{0\})$   \\
           & \hspace*{2pt}  1: $Cay(\Gamma, \{0\}\times \{ \pm 1\})$      &  \\
 \hline

  $(19,1)$ &  14: $P_1$   $(+1 \pmod {14}, -)$             & \hspace*{2pt} 1: $Cay(\Gamma, \{\pm 4\}\times \{0\})$ \\
             & \hspace*{2pt} 5:  Lemma~\ref{2p-533} with $a=1$ and $i=1$  &   \\
 \hline

\end{tabular}
\end{center}
}

Here are five methods to get $C_3$-factors. (1) From
$C_3$-factorization of certain Cayley graphs; (2) From several
initial $C_3$-factors $P_i$s by $(+7\pmod {14}, +1\pmod 3)$ or
$(+1\pmod {14}, -)$; (3) From several cycle sets $Q_i$s each of
which can generate a $C_3$-factor by $(+2\pmod {14}, -)$ since the
two elements having the same subscript in $Q_i$ have different
parity; (4) From  a cycle set $T$, each cycle of which can
generate a $C_3$-factor $F$ by $(+1 \pmod {14}, -)$, then 3
$C_3$-factors can be generated from $F$ by $(-,+1 \pmod {3})$; (5)
From a partial $C_3$-factor $S$, each cycle of $S$ will generate a
$C_3$-factor by $(+1 \pmod {14}, -)$.

For $C_7$-factors, we have the following four methods. (1) From
$C_7$-factorization of certain Cayley graphs; (2) From several
initial $C_7$-factors $P_i'$s by $(+7\pmod {14}, +1\pmod 3)$ or
$(+2\pmod {14}, -)$; (3) From a cycle set $Q'$, since these 7
elements having the same subscript in $Q'$ are different modulo 7,
so $Q'$ can generate a $C_7$-factor $F$ by $(+7 \pmod {14}, -)$,
then 3 $C_7$-factors can be generated from $F$ by $(-,+1 \pmod
{3})$; (4) From a cycle set $T'$, each cycle of which will
generate a $C_7$-factor by $(+7 \pmod {14}, +1 \pmod {3})$.

The cycles of $P_i$, $P_i'$, $Q_i$, $Q'$, $T$, $T'$ and $S$ are
given in  Appendix D. \qed

\begin{lemma} \label{RGDD84}
For any $(\alpha,\beta)\in \{(0,31),(9,22),(18,13),(24,7)\}$,
$(\alpha,\beta)\in {\rm HWP}(K_4[21];3,7)$.
\end{lemma}

\noindent {\it Proof:} Let the vertex set be $ Z_{21}\times Z_4$,
and the four parts of $K_4[21]$ be $Z_{21}\times \{i\}$, $i\in
Z_4$. For any $\alpha>0$, the required $\alpha$ $C_3$-factors  are
$$F_i^k=\{Q_i + (3j+k)_{0}\ |\  j=0,1,\cdots,6\},\ 1\le i\le
\alpha/3,\ k=0,1,2.$$ For $C_7$-factors, some of them will be
obtained from several cycle sets $Q_i'$. Here $\{Q_i' +
(7j+k)_{0}| j=0,1,2\}$
 is a $C_7$-factor for any $k=0,1,\cdots,6$
since these 7 elements having the same subscript in $Q_i'$ are
different modulo 7. The other $C_7$-factors will be obtained from
a cycle set $T'$ by $(+7 \pmod {21}, +1 \pmod {4})$, since the
first coordinate of the 7 elements in each cycle from $T'$ are
different modulo 7. For the sake of brevity, we list the 1-factor
$I$ and the cycles of $Q_i$, $Q_i'$ and $T'$ in Appendix E.\qed

\begin{lemma}
\label{L37-84} $(\alpha,\beta)\in {\rm HWP}(84;3,7)$ if and only
if $\alpha+\beta=41$.
\end{lemma}

\noindent {\it Proof:} Applying Construction~\ref{C-RGDD} with an
HW$(K_4[21]; 3, 7, \alpha_1, 31-\alpha_1)$,
 $\alpha_1=0,9,18,24$, from Lemma~\ref{RGDD84} and  an HW$(21; 3, 7, \alpha_2, 10-\alpha_2)$,
  $0\leq \alpha_2 \leq10$, from Theorem~\ref{L37} and Lemma~\ref{E21-37},
we can get an HW$( 84; 3, 7, \alpha_1 + \alpha_2,
41-\alpha_1-\alpha_2)$. Thus, $(\alpha, \beta)\in$ HWP$(84; 3, 7)$
for $0\le \alpha\le 34$. Similarly, for $(\alpha,\beta)=(35,6)$,
we can obtain the conclusion with an HW$(14; 3, 7, 0, 6)$ and an
HW$( K_6[14] ; 3, 7,  35,  0)$ from Theorem~\ref{Kv}. For all the
other cases, let the vertex set be $Z_{21}\times Z_4$. The methods
of generating the required $\alpha$ $C_3$-factors and $\beta$
$C_7$-factors are listed in Table 4.

For generating $C_3$-factors, here are four methods. (1) From a
$C_3$-factorization of certain Cayley graphs; (2) From an initial
$C_3$-factor $P$ by $(+ 1\pmod {21}, -)$; (3) From several cycle
sets $Q_i$s, note that $\{Q_i + (3j+k)_{0}| j=0,1,\cdots,6\}$
 is a $C_3$-factor for any $k=0,1,2$
since these 3 elements having the same subscript in $Q_i$ are
different modulo 3; (4) From a cycle set $T$ by $(+3 \pmod {21},
+1 \pmod {4})$, where the first coordinate of the 3 elements in
each cycle from $T$ are different modulo 3, so each cycle of $T$
will generate a $C_3$-factor by $(+3 \pmod {21}, +1 \pmod {4})$.
The required $C_7$-factors  are given from a $C_7$-factorization
of certain Cayley graphs or from a cycle set $T'$ by $(+7 \pmod
{21}, +1 \pmod {4})$. The cycles of $P$, $Q_i$, $T$ and $T'$ are
given in Appendix F. \qed

{\scriptsize
\begin{center}

\begin{tabular}{|c|l|l|l|}
\multicolumn{4}{c}{\bf Table 4 HWP$(84;3,7)$}\\[5pt]
  \hline
 $(\alpha, \beta)$  &  \hspace*{1cm} $ C_3$-factor &  $  \hspace*{1cm}  C_7$-factor  &  \hspace*{1.5cm} 1-factor \\   \hline

   $(36,5)$ & 21: $P$     &  4: $T'$                                   & $\{(i_0,(i+17)_1),(i_2,(i+4)_3)| i \in Z_{21}\}$  \\
            & 15: $Q_i$, $1\leq i \leq 5 $        &  1: $Cay(\Gamma, \{\pm 9\}\times \{0\})$  & \\
 \hline

           &  21: $P$                &   &\\
  $(37,4)$ &  15: $Q_i$, $1\leq i \leq 5 $                  &  4: $T'$ & $\{(i_0,(i+16)_1),(i_2,(i+12)_3)| i \in Z_{21}\}$ \\
           &  \hspace*{2pt} 1: $Cay(\Gamma, \{\pm 7\}\times \{0\})$      &   & \\
 \hline

             & 21: $P$   &  1: $Cay(\Gamma, \{\pm 3\}\times \{0\})$  &  \\
  $(38,3)$   & 12: $Q_i$, $1\leq i \leq 4 $          &  1: $Cay(\Gamma, \{\pm 6\}\times \{0\})$  & $\{(i_0,(i+1)_1),(i_2,(i+13)_3)| i \in Z_{21}\}$ \\
             & \hspace*{2pt} 5: $T$                         &  1: $Cay(\Gamma, \{\pm 9\}\times \{0\})$  &  \\
 \hline

  $(39,2)$  & 21: $P$      &  1: $Cay(\Gamma, \{\pm 6\}\times \{0\})$  & $\{(i_0,(i+20)_1),(i_2, i_3)| i \in Z_{21}\}$  \\
            & 18: $Q_i$, $1\leq i \leq 6 $        &  1: $Cay(\Gamma, \{\pm 9\}\times \{0\})$  & \\
 \hline

           &  21: $P$              &   &\\
  $(40,1)$ &  18: $Q_i$, $1\leq i \leq 6 $                  &  1: $Cay(\Gamma, \{\pm 9\}\times \{0\})$  & $\{(i_0,(i+5)_1),(i_2,(i+11)_3)| i \in Z_{21}\}$ \\
           &  \hspace*{2pt} 1: $Cay(\Gamma, \{\pm 7\}\times \{0\})$      &   & \\
 \hline

\end{tabular}
\end{center}
}

\begin{lemma}
\label{L37even}  If $v \equiv 0\pmod {42}$, then
$(\alpha,\beta)\in {\rm HWP}(v;3,7)$ with
$\alpha+\beta=\frac{v-2}{2}$.
\end{lemma}

\noindent {\it Proof:}
Let $v=42u$, $u\ge 1$.
For $u\le 2$, the
conclusion comes from Lemmas~\ref{L37-42} and ~\ref{L37-84}.

For $u=3$, start with an HW$(K_3[3];3,7,3,0)$, an
HW$(C_3[14];3,7,14,0)$ and an HW$(C_3[14];$ $3,7,0,14)$ from
Theorem~\ref{Kv}, apply Construction~\ref{L351} with $s=14$ and
$t_i\in\{0,14\}$ to get an HW$(K_3[42];$ $
3,7,\sum_{i=1}^{3}t_i,42-\sum_{i=1}^{3}t_i)$. Further, applying
Construction~\ref{C-RGDD} with an HW$(42;3,7,\alpha',20-\alpha')$,
$0\leq \alpha'\leq20$, from Lemma~\ref{L37-42} to obtain an
HW$(126;3,7,\sum_{i=1}^{3}t_i+\alpha',42-\sum_{i=1}^{3}t_i+(20-\alpha'))$.
Thus we have obtained an HW$(126;3,7,\alpha,\beta)$ for any
$\alpha+\beta=62$ since $\sum_{i=1}^{3}t_i+\alpha'$ can cover the
integers from $0$ to $62$.

For $u\geq 4$, similarly, start with an HW$(K_u[6];3,7,3u-3,0)$,
an HW$(C_3[7];3,7,7,0)$ and an HW$(C_3[7];3,7,0,7)$ from
Theorem~\ref{Kv}, and apply Construction~\ref{L351} with $s=7$ and
$t_i\in\{0,7\}$ to get an
HW$(K_u[42];3,7,\sum_{i=1}^{3u-3}t_i,21u-21-\sum_{i=1}^{3u-3}t_i)$.
Further, applying Construction~\ref{C-RGDD} with an
HW$(42;3,7,\alpha',20-\alpha')$, $0\leq \alpha'\leq20$, from
Lemma~\ref{L37-42}, we can obtain an
HW$(42u;3,7,\alpha'+\sum_{i=1}^{3u-3}t_i,20-\alpha'+21u-21-\sum_{i=1}^{3u-3}t_i)$.
It's easy to prove that $\alpha'+\sum_{i=1}^{3u-3}t_i$ can cover
the integers from $0$ to $21u-1$. The proof is complete.\qed

Combining Theorem~\ref{L37}, Lemmas~\ref{E21-37} and ~\ref{L37even},
we have the following theorem.

\begin{theorem}
\label{TH37} $ (\alpha, \beta)\in {\rm HWP}(v;3,7)$ if and only if
$v\equiv 0\pmod{21}$ and $\alpha+
\beta=\lfloor\frac{v-1}{2}\rfloor$.
\end{theorem}

Combining Theorems~\ref{TH34}, ~\ref{TH35} and ~\ref{TH37}, we
have proved Theorem~\ref{TH457}.

\newpage

\section*{Appendix A for Lemma~\ref{L35-30} }

{\footnotesize  \noindent
 {$(\alpha,\beta)=(1,13) :$}\\[5pt]
 \noindent
\begin{tabular}{lllllllllll}
$P_1'$
&$( 0_1, 2_0, 4_2, 1_2, 3_1)$
&$( 5_2, 6_0, 7_1, 8_2, 9_0)$
&$( 5_0, 7_2, 9_1, 6_1, 8_0)$\\
&$( 0_0, 1_1, 2_2, 3_0, 5_1)$
&$( 4_1, 0_2, 6_2, 1_0, 7_0)$
&$( 2_1, 8_1, 3_2, 4_0, 9_2)$\\
$P_2'$
&$( 0_0, 4_1, 5_2, 1_1, 8_2)$
&$( 2_2, 6_0, 3_0, 7_1, 0_1)$
&$( 9_0, 1_2, 9_1, 3_1, 0_2)$\\
&$( 2_0, 5_0, 1_0, 4_2, 8_1)$
&$( 6_1, 3_2, 5_1, 8_0, 9_2)$
&$( 7_2, 2_1, 7_0, 4_0, 6_2)$\\
$P_3'$
&$( 0_0, 5_2, 2_0, 1_1, 7_2)$
&$( 2_2, 1_2, 3_0, 4_2, 3_2)$
&$( 4_1, 5_0, 9_2, 8_2, 5_1)$\\
&$( 6_0, 2_1, 9_1, 4_0, 8_1)$
&$( 7_1, 8_0, 1_0, 0_1, 6_2)$
&$( 9_0, 3_1, 7_0, 0_2, 6_1)$\\
$P_4'$
&$( 0_0, 5_2, 2_0, 1_1, 7_2)$
&$( 2_2, 1_2, 3_0, 4_2, 3_2)$
&$( 4_1, 5_0, 9_2, 8_2, 5_1)$\\
&$( 6_0, 2_1, 9_1, 4_0, 8_1)$
&$( 7_1, 8_0, 1_0, 0_1, 6_2)$
&$( 9_0, 3_1, 7_0, 0_2, 6_1)$
\end{tabular}}

{\footnotesize \vspace{5pt}\noindent $(\alpha,\beta)=(2,12):$\\[5pt]
\noindent
\begin{tabular}{lllllllllll}
$Q_1$ &$( 0_0, 1_1, 2_2)$ &$( 3_0, 4_1, 5_2)$& $Q_2$ &$( 0_0, 3_0,
7_1)$ &$( 2_2, 4_1, 1_2)$
\end{tabular}}

{\footnotesize  \noindent
\begin{tabular}{lllllllll}
$P_1'$
&$( 0_0, 4_1, 1_1, 3_0, 8_2)$
&$( 2_2, 5_2, 6_0, 1_2, 7_1)$
&$( 9_0, 2_0, 0_1, 3_1, 8_0)$\\
&$( 4_2, 5_0, 8_1, 0_2, 4_0)$
&$( 6_1, 1_0, 9_2, 9_1, 7_0)$
&$( 7_2, 2_1, 6_2, 3_2, 5_1)$\\
$P_2'$
&$( 0_0, 9_0, 5_2, 4_0, 7_2)$
&$( 1_1, 8_2, 7_0, 1_2, 0_2)$
&$( 2_2, 5_0, 2_1, 6_0, 9_1)$\\
&$( 3_0, 3_1, 9_2, 0_1, 3_2)$
&$( 4_1, 2_0, 6_2, 6_1, 5_1)$
&$( 7_1, 8_0, 8_1, 4_2, 1_0)$
\end{tabular}}

{\footnotesize \vspace{5pt}
\noindent {$(\alpha,\beta)=(3,11) : $}\\[5pt]
\noindent
\begin{tabular}{lllllllllll}
$Q_1$ &$( 0_0, 1_1, 2_2)$ &$( 3_0, 4_1, 5_2)$ & $Q_2$ &$( 0_0,3_0,
7_1)$ &$( 2_2, 4_1, 1_2)$
\end{tabular}}

{\footnotesize  \noindent
\begin{tabular}{lllllllll}
$P_1'$
&$( 0_0, 4_1, 1_1, 3_0, 6_0)$
&$( 2_2, 5_2, 7_1, 1_2, 5_0)$
&$( 8_2, 9_0, 4_2, 0_1, 2_0)$\\
&$( 3_1, 8_0, 6_2, 7_2, 0_2)$
&$( 6_1, 2_1, 5_1, 9_1, 7_0)$
&$( 1_0, 8_1, 3_2, 4_0, 9_2)$\\
$P_2'$
&$( 0_0, 9_0, 2_2, 8_2, 2_1)$
&$( 1_1, 5_0, 4_0, 1_2, 8_0)$
&$( 3_0, 4_2, 9_1, 0_1, 7_0)$\\
&$( 4_1, 3_2, 7_2, 6_0, 5_1)$
&$( 5_2, 3_1, 9_2, 6_1, 1_0)$
&$( 7_1, 0_2, 8_1, 2_0, 6_2)$
\end{tabular}}

{\footnotesize \vspace{5pt}
\noindent {$(\alpha,\beta)=(4,10) : $}\\[5pt]
\noindent
\begin{tabular}{lllllllllll}
$Q_1$ &$( 0_0, 1_1, 2_2)$ &$( 3_0, 4_1, 5_2)$ &$Q_2$ &$( 0_0, 3_0,
7_1)$
&$( 2_2, 4_1, 1_2)$\\
$Q_3$ &$( 0_0, 4_1, 8_2)$ &$( 1_1, 3_0, 1_2)$& $Q_4$ &$( 0_0, 5_2,
9_0)$ &$( 1_1, 4_1, 4_2)$
\end{tabular}}

{\footnotesize  \noindent
\begin{tabular}{lllllllll}
$P_1'$
&$( 0_0, 6_0, 1_1, 5_2, 0_1)$
&$( 2_2, 3_0, 8_2, 1_2, 2_0)$
&$( 4_1, 7_1, 3_1, 0_2, 9_0)$\\
&$( 4_2, 8_0, 9_2, 3_2, 6_2)$
&$( 5_0, 1_0, 4_0, 2_1, 8_1)$
&$( 6_1, 5_1, 7_2, 9_1, 7_0)$\\
$P_2'$
&$( 0_0, 2_0, 1_1, 0_2, 4_2)$
&$( 2_2, 5_0, 7_0, 6_0, 3_2)$
&$( 3_0, 0_1, 8_0, 4_1, 6_2)$\\
&$( 5_2, 7_2, 4_0, 7_1, 1_0)$
&$( 8_2, 3_1, 5_1, 1_2, 2_1)$
&$( 9_0, 9_1, 8_1, 6_1, 9_2)$
\end{tabular}}

{\footnotesize \vspace{5pt}
\noindent {$(\alpha,\beta)=(5,9) : $}\\[5pt]
\noindent
\begin{tabular}{lllllllllll}
$P_1'$
&$( 0_0, 3_0, 6_0, 1_1, 4_1)$
&$( 2_2, 5_2, 8_2, 1_2, 7_1)$
&$( 9_0, 2_0, 6_1, 0_1, 5_0)$\\
&$( 3_1, 7_2, 0_2, 4_2, 6_2)$
&$( 8_0, 2_1, 8_1, 9_1, 3_2)$
&$( 1_0, 5_1, 9_2, 4_0, 7_0)$\\
$P_2'$
&$( 0_0, 7_1, 2_0, 3_0, 9_0)$
&$( 1_1, 8_2, 5_0, 0_2, 3_1)$
&$( 2_2, 9_1, 4_0, 7_2, 6_2)$\\
&$( 4_1, 1_2, 5_1, 6_1, 1_0)$
&$( 5_2, 2_1, 9_2, 6_0, 8_0)$
&$( 0_1, 3_2, 8_1, 4_2, 7_0)$\\
$P_3'$
&$( 0_0, 2_0, 5_2, 3_2, 4_2)$
&$( 1_1, 0_1, 6_2, 3_0, 9_1)$
&$( 2_2, 1_2, 8_1, 4_0, 8_0)$\\
&$( 4_1, 2_1, 6_0, 9_0, 5_1)$
&$( 7_1, 3_1, 7_0, 5_0, 1_0)$
&$( 8_2, 7_2, 9_2, 6_1, 0_2)$
\end{tabular}}

{\footnotesize \vspace{5pt}\noindent {$(\alpha,\beta)=(6,8) : $}\\[5pt]
\noindent
\begin{tabular}{lllllllllll}
$P_1$
&$( 0_2, 4_0, 7_0)$
&$( 3_0, 4_1, 5_2)$
&$( 6_0, 7_1, 8_2)$
&$( 9_0, 0_1, 3_1)$
&$( 1_2, 4_2, 6_1)$\\
&$( 2_0, 5_0, 7_2)$
&$( 8_0, 9_1, 1_0)$
&$( 2_1, 6_2, 9_2)$
&$( 0_0, 1_1, 3_2)$
&$( 2_2, 5_1, 8_1)$\\
$P_2$
&$( 0_0, 2_2, 0_1)$
&$( 1_1, 3_0, 8_2)$
&$( 4_1, 9_0, 2_0)$
&$( 5_2, 6_0, 3_1)$
&$( 7_1, 1_2, 0_2)$\\
&$( 4_2, 1_0, 5_1)$ &$( 5_0, 9_1, 8_1)$ &$( 6_1, 7_0, 9_2)$ &$(
7_2, 2_1, 4_0)$ &$( 8_0, 3_2, 6_2)$
\end{tabular}}

{\footnotesize \vspace{2pt}\noindent
\begin{tabular}{lllllllll}
$P_1'$
&$( 0_0, 5_2, 0_1, 6_0, 1_2)$
&$( 1_1, 2_2, 3_0, 7_1, 5_0)$
&$( 4_1, 7_2, 7_0, 8_0, 5_1)$\\
&$( 8_2, 9_1, 6_2, 6_1, 4_0)$
&$( 9_0, 4_2, 3_2, 3_1, 2_1)$
&$( 2_0, 1_0, 9_2, 0_2, 8_1)$\\
$P_2'$
&$( 0_0, 3_1, 6_2, 7_2, 8_1)$
&$( 1_1, 4_2, 8_0, 8_2, 1_0)$
&$( 2_2, 2_0, 5_1, 5_2, 6_1)$\\
&$( 3_0, 9_1, 7_0, 1_2, 2_1)$
&$( 4_1, 3_2, 7_1, 0_1, 4_0)$
&$( 6_0, 5_0, 9_2, 9_0, 0_2)$
\end{tabular}}

\newpage
{\footnotesize \vspace{5pt}
\noindent {$(\alpha,\beta)=(7,7) : $}\\[5pt]
\noindent
\begin{tabular}{lllllllllll}
$P_1$
&$( 0_2, 4_0, 7_0)$
&$( 3_0, 4_1, 5_2)$
&$( 6_0, 7_1, 8_2)$
&$( 9_0, 0_1, 3_1)$
&$( 1_2, 4_2, 6_1)$\\
&$( 2_0, 5_0, 7_2)$
&$( 8_0, 9_1, 1_0)$
&$( 0_0, 1_1, 2_2)$
&$( 2_1, 3_2, 8_1)$
&$( 5_1, 6_2, 9_2)$\\
$P_2$
&$( 0_0, 5_2, 8_2)$
&$( 1_1, 3_0, 7_1)$
&$( 2_2, 4_1, 9_0)$
&$( 6_0, 0_1, 9_1)$
&$( 1_2, 8_0, 0_2)$\\
&$( 2_0, 3_2, 9_2)$
&$( 3_1, 6_1, 7_0)$
&$( 4_2, 1_0, 5_1)$
&$( 5_0, 2_1, 6_2)$
&$( 7_2, 4_0, 8_1)$\\
\end{tabular}}

{\footnotesize \noindent
\begin{tabular}{lllllllll}
$P_1'$
&$( 0_0, 6_0, 1_2, 2_2, 6_1)$
&$( 1_1, 9_0, 2_0, 3_0, 5_1)$
&$( 4_1, 0_1, 4_0, 8_0, 3_2)$\\
&$( 5_2, 6_2, 9_1, 7_0, 9_2)$
&$( 7_1, 3_1, 0_2, 7_2, 1_0)$
&$( 8_2, 5_0, 8_1, 4_2, 2_1)$\\
$P_2'$
&$( 0_0, 1_2, 4_0, 3_0, 9_1)$
&$( 1_1, 2_0, 0_2, 7_1, 5_0)$
&$( 2_2, 7_0, 4_2, 6_0, 8_1)$\\
&$( 4_1, 5_1, 9_0, 8_0, 9_2)$
&$( 5_2, 0_1, 3_2, 6_1, 2_1)$
&$( 8_2, 7_2, 6_2, 3_1, 1_0)$
\end{tabular}}

{\footnotesize \vspace{5pt}
\noindent {$(\alpha,\beta)=(8,6) : $}\\
\noindent
\begin{tabular}{lllllllllll}
$Q_1$ &$( 0_0, 1_1, 2_2)$ &$( 3_0, 4_1, 5_2)$ &$Q_2$ &$( 0_0, 3_0,
7_1)$ &$( 2_2, 4_1, 1_2)$ & $Q_3$ &$( 0_0, 4_1, 8_2)$ &$( 1_1,
3_0, 1_2)$
\end{tabular}}

{\footnotesize  \noindent
\begin{tabular}{lllllllllll}
$P_1$ &$( 0_0, 5_2, 6_0)$ &$( 1_1, 4_1, 7_1)$ &$( 2_2, 3_0, 8_2)$
&$( 9_0, 2_0, 2_1)$
&$( 0_1, 5_0, 1_0)$\\
&$( 1_2, 3_1, 4_0)$ &$( 4_2, 8_0, 7_0)$ &$( 6_1, 6_2, 9_2)$ &$(
7_2, 0_2, 5_1)$ &$( 9_1, 3_2, 8_1)$
\end{tabular}}

{\footnotesize \noindent
\begin{tabular}{lllllllllll}
$P_1'$
&$( 0_0, 1_2, 5_2, 2_0, 5_1)$
&$( 1_1, 5_0, 7_1, 7_0, 0_2)$
&$( 2_2, 1_0, 7_2, 3_0, 3_2)$\\
&$( 4_1, 6_0, 2_1, 4_2, 8_1)$
&$( 8_2, 8_0, 9_0, 6_1, 4_0)$
&$( 0_1, 9_1, 6_2, 3_1, 9_2)$
\end{tabular}}

{\footnotesize \vspace{5pt}
\noindent {$(\alpha,\beta)=(9,5) : $}\\
\noindent
\begin{tabular}{lllllllllll}
$Q_1$ &$( 0_0, 1_1, 2_2)$ &$( 3_0, 4_1, 5_2)$ &$Q_2$ &$( 0_0, 3_0,
7_1)$ &$( 2_2, 4_1, 1_2)$ & $Q_3$ &$( 0_0, 4_1, 8_2)$
&$( 1_1, 3_0, 7_2)$\\
\end{tabular}}

{\footnotesize  \noindent
\begin{tabular}{lllllllllll}
$P_1$ &$( 0_0, 5_2, 6_0)$ &$( 1_1, 4_1, 7_1)$ &$( 2_2, 3_0, 8_2)$
&$( 9_0, 2_0, 1_0)$
&$( 0_1, 5_0, 4_0)$\\
&$( 1_2, 3_1, 8_0)$ &$( 4_2, 9_1, 6_2)$ &$( 6_1, 5_1, 9_2)$ &$(
7_2, 2_1, 8_1)$ &$( 0_2, 3_2, 7_0)$
\end{tabular}}

{\footnotesize \noindent
\begin{tabular}{lllllllll}
$P_1'$
&$( 0_0, 1_2, 3_0, 6_1, 8_1)$
&$( 1_1, 8_0, 0_1, 2_2, 9_1)$
&$( 4_1, 5_0, 7_1, 1_0, 7_0)$\\
&$( 5_2, 4_2, 5_1, 7_2, 9_2)$
&$( 6_0, 0_2, 9_0, 6_2, 3_2)$
&$( 8_2, 2_0, 4_0, 3_1, 2_1)$
\end{tabular}}

{\footnotesize \vspace{5pt}
\noindent {$(\alpha,\beta)=(11,3) : $}\\[5pt]
\noindent
\begin{tabular}{lllllllllll}
$P_1$
&$( 0_0, 3_0, 6_0)$
&$( 1_1, 4_1, 7_1)$
&$( 2_2, 5_2, 8_2)$
&$( 9_0, 2_0, 6_1)$
&$( 0_1, 4_2, 9_1)$\\
&$( 1_2, 0_2, 5_1)$
&$( 3_1, 2_1, 9_2)$
&$( 5_0, 4_0, 8_1)$
&$( 7_2, 1_0, 6_2)$
&$( 8_0, 3_2, 7_0)$\\
$P_2$
&$( 0_0, 7_1, 3_1)$
&$( 1_1, 8_0, 5_1)$
&$( 2_2, 5_0, 7_0)$
&$( 3_0, 6_1, 1_0)$
&$( 4_1, 1_2, 9_2)$\\
&$( 5_2, 9_0, 6_2)$ &$( 6_0, 0_1, 2_1)$ &$( 8_2, 4_2, 0_2)$ &$(
2_0, 7_2, 4_0)$ &$( 9_1, 3_2, 8_1)$
\end{tabular}}

{\footnotesize \noindent
\begin{tabular}{lllllllll}
$P_1'$
&$( 0_0, 4_2, 6_2, 2_2, 5_1)$
&$( 1_1, 5_2, 9_1, 6_1, 8_1)$
&$( 3_0, 0_1, 7_0, 9_0, 5_0)$\\
&$( 4_1, 3_1, 6_0, 0_2, 7_2)$
&$( 7_1, 1_0, 2_0, 9_2, 3_2)$
&$( 8_2, 2_1, 8_0, 1_2, 4_0)$
\end{tabular}}

{\footnotesize \vspace{5pt}
\noindent {$(\alpha,\beta)=(12,2) : $}\\[5pt]
\noindent
\begin{tabular}{lllllllllll}
$P_1$
&$( 0_0, 1_1, 2_2)$
&$( 3_0, 4_1, 5_2)$
&$( 6_0, 9_0, 3_1)$
&$( 7_1, 0_1, 2_0)$
&$( 8_2, 1_2, 5_0)$\\
&$( 4_2, 6_1, 4_0)$
&$( 7_2, 5_1, 7_0)$
&$( 8_0, 2_1, 9_2)$
&$( 9_1, 3_2, 6_2)$
&$( 0_2, 1_0, 8_1)$\\
$P_2$
&$( 0_0, 5_2, 0_1)$
&$( 1_1, 7_2, 8_0)$
&$( 2_2, 5_0, 4_0)$
&$( 3_0, 1_2, 3_1)$
&$( 4_1, 4_2, 3_2)$\\
&$( 6_0, 0_2, 5_1)$ &$( 7_1, 1_0, 6_2)$ &$( 8_2, 2_1, 7_0)$ &$(
9_0, 2_0, 8_1)$ &$( 6_1, 9_1, 9_2)$
\end{tabular}}

{\footnotesize \noindent
\begin{tabular}{lllllllllll}
 $Q_1$ &$( 0_0, 9_0, 3_2)$ &$(
1_1, 0_1, 4_2)$ &$Q_2$ &$( 0_0, 7_2, 6_2)$ &$( 1_1, 9_0, 2_1)$
\end{tabular}}

{\footnotesize
\noindent {$(\alpha,\beta)=(13,1) : $}\\
\noindent
\begin{tabular}{lllllllllll}
$Q_1$ &$( 0_0, 8_2, 1_0)$ &$( 1_1, 0_1, 3_2)$& $Q_2$ &$( 0_0, 0_1,
0_2)$ &$( 1_1, 7_2, 1_0)$ & $Q_3$ &$( 0_0, 1_2, 3_1)$ &$( 2_2,
0_1, 1_0)$
\end{tabular}}

{\footnotesize  \noindent
\begin{tabular}{lllllllllll}
$P_1$

&$( 0_0, 1_1, 2_2)$ &$( 3_0, 4_1, 5_2)$ &$( 6_0, 9_0, 2_0)$ &$(
7_1, 0_1, 3_1)$
&$( 8_2, 1_2, 4_2)$\\
&$( 5_0, 9_1, 1_0)$ &$( 6_1, 8_0, 3_2)$ &$( 7_2, 4_0, 8_1)$ &$(
0_2, 5_1, 9_2)$
&$( 2_1, 6_2, 7_0)$\\
$P_2$ &$( 0_0, 7_1, 4_2)$ &$( 1_1, 0_2, 2_1)$ &$( 2_2, 9_0, 6_1)$
&$( 3_0, 9_1, 4_0)$
&$( 4_1, 0_1, 8_0)$\\
&$( 5_2, 6_0, 5_1)$ &$( 8_2, 2_0, 9_2)$ &$( 1_2, 7_2, 7_0)$ &$(
3_1, 1_0, 6_2)$ &$( 5_0, 3_2, 8_1)$

\end{tabular}}

\section*{Appendix B for Lemma~\ref{RGDD60}}

{\footnotesize \vspace{5pt} \noindent {$(\alpha,\beta)=(0,22) : $
\ \ $I=\{(i_0,(i+13)_3),(i_1,(i+2)_2)| i \in Z_{15}\}$. }\\[5pt]
\noindent
\begin{tabular}{lllllllllll}
$Q_1'$
&$( 0_0, 1_1, 2_2, 3_3, 5_1)$
&$( 4_0, 6_2, 8_0, 7_3,14_1)$
&$(10_2, 0_3, 6_0, 6_3, 2_0)$
&$(13_1, 3_2, 7_1,14_3, 4_2)$\\
$Q_2'$ &$( 0_0,10_2,12_3,10_1,10_3)$ &$( 1_1,14_0, 6_2, 1_0, 6_3)$
&$( 2_2, 8_3,13_1, 7_0, 7_1)$ &$( 8_0, 8_2, 4_3, 4_2, 4_1)$
\end{tabular}}

{\footnotesize\noindent
\begin{tabular}{lllllllllll}
$T'$ &$( 0_0,12_3, 9_1, 8_2, 6_3)$ &$( 0_0, 8_3,11_1, 9_2, 2_3)$
&$( 0_0, 6_2, 3_1, 7_3, 9_2)$
&$( 0_0, 9_1, 6_0,13_1,12_2)$\\
&$( 0_0,14_2, 7_1,13_3, 1_2)$
&$( 0_0, 7_3, 1_0, 9_3, 8_1)$
&$( 0_0,11_2,12_1, 8_0, 4_2)$
&$( 0_0, 3_2, 9_0,12_1, 1_3)$\\
&$( 0_0, 3_1,14_0,11_2, 7_1)$ &$( 0_0,13_1, 7_2, 6_0,14_1)$ &$(
0_0, 4_1,13_2, 6_3,12_1)$ &$( 0_0, 3_3, 2_1,14_2,11_3)$
\end{tabular}}

{\footnotesize \vspace{5pt}\noindent {$(\alpha,\beta)=(6,16) : $\
\ $I=\{(i_0,(i+4)_2),(i_1,(i+4)_3)| i \in Z_{15}\}$.}\\
\noindent
\begin{tabular}{lllllllllll}
$T$ &$( 0_0,14_2,10_1)$ &$( 0_0,10_2, 2_3)$ &$( 0_0, 5_3, 4_1)$
&$( 0_0,11_3, 1_2)$ &$( 0_0,13_1, 8_3)$ &$( 0_0, 7_1, 5_2)$
\end{tabular}}

{\footnotesize \noindent
\begin{tabular}{lllllllllll}
$Q_1'$
&$( 0_0, 5_1,10_2, 3_3, 3_2)$
&$( 1_1, 8_0, 0_3, 2_0,12_3)$
&$( 2_2, 9_1, 9_0, 6_2, 3_1)$
&$(11_3, 9_2,14_3, 2_1,11_0)$\\
$Q_2'$
&$( 0_0, 2_1, 4_2,11_0,11_2)$
&$( 1_1, 7_2,13_3, 5_1, 2_0)$
&$( 4_0, 1_3,13_2,13_1, 5_3)$
&$( 7_3,13_0, 5_2, 4_3, 4_1)$\\
$T'$
&$( 0_0, 6_2,12_1, 8_2, 4_3)$
&$( 0_0,12_1,11_0,14_3,13_2)$
&$( 0_0, 6_3, 8_1, 2_2,14_3)$
&$( 0_0,11_1, 2_3, 3_0, 9_2)$\\
&$( 0_0, 2_2, 6_1, 3_2, 9_1)$
&$( 0_0, 1_1,14_3,12_1, 3_3)$
\end{tabular}}

{\footnotesize \vspace{5pt}
\noindent {$(\alpha,\beta)=(12,10) : $\ \ $I=\{(i_0,(i+14)_1),(i_2,(i+3)_3)| i \in Z_{15}\}$.}\\[5pt]
\noindent
\begin{tabular}{lllllllllll}
$T$
&$( 0_0, 1_1, 2_2)$
&$( 0_0, 2_3,13_2)$
&$( 0_0, 7_1, 5_2)$
&$( 0_0, 8_2,13_3)$
&$( 0_0,11_3, 4_1)$
&$( 0_0, 5_1, 1_2)$\\
&$( 0_0, 7_3,14_2)$ &$( 0_0, 7_2,14_3)$ &$( 0_0,11_1,10_3)$ &$(
0_0, 2_1, 4_3)$ &$( 0_0,13_1, 8_3)$ &$( 0_0,10_2, 8_1)$
\end{tabular}}

{\footnotesize \noindent
\begin{tabular}{lllllllllll}
$Q_1'$
&$( 0_0, 3_3, 6_2, 9_1, 0_3)$
&$( 1_1, 4_0,10_2, 1_0, 4_2)$
&$( 2_2,11_3, 2_1, 8_0, 2_3)$
&$( 5_1, 2_0, 8_1, 3_2, 9_3)$\\
$Q_2'$
&$( 0_0,10_1, 4_2,13_1,12_2)$
&$( 1_1,12_3,12_1, 0_3, 1_2)$
&$( 3_3,12_0, 9_3, 8_0, 8_2)$
&$( 4_0, 0_2,11_0, 1_3, 4_1)$
\end{tabular}}

\section*{Appendix C for Lemma~\ref{L35-60}}

{\footnotesize
\noindent {$(\alpha,\beta)=(18,11) : $}\\[5pt]
\noindent
\begin{tabular}{lllllllllll}
$P$
&$(14_2,13_3, 7_0)$
&$( 0_0,14_0,10_3)$
&$( 2_2, 4_3,11_1)$
&$( 8_0,10_1,12_1)$
&$( 3_3,14_1, 3_0)$\\
&$( 6_2, 8_2, 8_1)$
&$(10_2, 2_0,13_2)$
&$( 9_1, 3_2, 4_2)$
&$( 5_0, 6_1, 9_3)$
&$( 5_1,12_3,10_0)$\\
&$(13_0,12_2,11_0)$
&$( 7_3, 0_3, 5_3)$
&$(12_0, 1_2, 6_3)$
&$( 4_0, 1_3, 2_3)$
&$(11_3, 7_2, 0_2)$\\
&$(13_1, 8_3, 3_1)$
&$( 2_1, 5_2,14_3)$
&$( 1_1, 9_0, 0_1)$
&$( 1_0,11_2, 6_0)$
&$( 7_1, 4_1, 9_2)$\\
$Q_1$ &$( 0_0, 8_0, 1_3)$ &$( 1_1,12_3, 7_0)$ &$( 2_2, 6_1, 5_3)$
&$( 5_1, 4_2, 9_2)$
\end{tabular}}

{\footnotesize \noindent
\begin{tabular}{lllllllll}
$Q'$
&$( 0_0, 5_1, 6_2, 4_3, 3_1)$
&$( 1_1, 2_0, 1_3, 0_2, 3_0)$
&$( 2_2, 8_3, 8_2, 0_3,12_3)$
&$( 4_0, 1_0, 4_2,12_1, 4_1)$\\
$T'$
&$( 0_0, 2_2, 8_0, 11_3, 4_0)$
&$( 0_0, 6_2, 9_1, 2_2, 13_2)$
&$( 0_0, 11_0, 3_3, 14_3,12_1)$
&$( 0_0, 7_3, 1_1,14_3, 8_1)$\\
&$( 0_0, 3_3, 1_1,12_1, 9_2)$
\end{tabular}}

{\footnotesize \vspace{5pt}
\noindent {$(\alpha,\beta)=(21,8) : $}\\[5pt]
\noindent
\begin{tabular}{lllllllllll}
$P$
&$(10_2, 3_1, 5_2)$
&$(13_1, 7_2, 9_0)$
&$( 1_0, 2_0,12_1)$
&$( 2_2,12_3, 1_2)$
&$(10_1,10_0, 6_3)$\\
&$(12_0, 3_2, 1_3)$
&$( 1_1, 5_3,13_3)$
&$(11_2, 6_0,13_2)$
&$( 5_0, 4_2,11_0)$
&$( 2_1, 4_3,12_2)$\\
&$( 0_3, 2_3, 8_1)$
&$( 5_1,11_1, 9_2)$
&$( 4_0,11_3,10_3)$
&$( 7_3, 9_1, 7_0)$
&$( 6_1, 7_1,14_0)$\\
&$( 6_2,14_2, 3_0)$
&$( 0_0, 8_0,14_1)$
&$(13_0, 8_2,14_3)$
&$( 3_3, 8_3, 0_1)$
&$( 0_2, 9_3, 4_1)$\\
$Q_1$
&$( 0_0, 2_2, 5_0)$
&$( 1_1,11_1, 2_3)$
&$( 3_3, 4_2,10_0)$
&$( 6_2, 0_1,10_3)$\\
$T$
&$( 0_0, 4_0, 5_1)$\\
\end{tabular}}

{\footnotesize \noindent
\begin{tabular}{lllllllll}
$Q'$
&$( 0_0, 9_1, 8_2, 8_3, 2_3)$
&$( 1_1, 3_1,10_1,12_0, 1_2)$
&$( 2_2, 4_3,14_0, 1_0,10_3)$
&$( 8_0, 6_3, 7_1, 0_2, 9_2)$\\
$T'$ &$( 0_0, 3_3, 6_3, 9_3, 12_3)$ &$( 0_0, 12_0, 9_1,6_2, 3_1)$
&$( 0_0, 3_0,6_3, 9_0,12_1)$

\end{tabular}}

{\footnotesize \vspace{5pt}
\noindent {$(\alpha,\beta)=(22,7) : $}\\[5pt]
\noindent
\begin{tabular}{lllllllllll}
$P$
&$( 2_1, 8_2,10_0)$
&$( 2_2, 9_0, 0_1)$
&$(13_1, 7_2,14_3)$
&$(14_2, 5_0, 1_2)$
&$( 5_1, 1_3,12_1)$\\
&$( 0_0,13_0,12_2)$
&$(11_3, 9_3,11_0)$
&$(12_3, 0_2, 4_1)$
&$( 6_1, 7_1, 5_2)$
&$( 8_3, 3_1, 7_0)$\\
&$( 3_3, 1_0,11_1)$
&$(14_1, 6_0,13_3)$
&$( 8_0, 4_3,13_2)$
&$( 7_3, 3_2,11_2)$
&$( 4_0, 8_1,10_3)$\\
&$(12_0, 4_2, 6_3)$
&$( 1_1, 2_0, 3_0)$
&$(10_2,10_1, 5_3)$
&$( 6_2, 9_1, 0_3)$
&$(14_0, 2_3, 9_2)$\\
$Q_1$
&$( 0_0, 2_2, 9_1)$
&$( 1_1, 1_0,13_3)$
&$( 3_3, 4_2, 2_3)$
&$( 5_1, 8_0,12_2)$\\
$T$
&$( 0_0, 4_0, 5_1)$\\
\end{tabular}}

{\footnotesize \noindent
\begin{tabular}{lllllllll}
$Q'$
&$( 0_0, 6_2, 3_1, 1_1, 4_3)$
&$( 2_2, 3_2, 2_0, 7_3,14_0)$
&$( 3_3,11_3, 2_1, 0_3, 0_2)$
&$( 5_1, 3_0,11_0,14_1, 9_2)$
\end{tabular}}

{\footnotesize \vspace{5pt}
\noindent {$(\alpha,\beta)=(23,6) : $}\\[5pt]
 \noindent
\begin{tabular}{lllllllll}
$Q_1$ &$( 0_0, 5_1, 9_1)$ &$( 1_1, 8_0,11_2)$ &$( 3_3,10_2,10_3)$
&$( 4_0, 0_2, 2_3)$\\
$S$ &$( 0_0, 1_1, 2_2)$ &$( 0_0, 4_3, 11_1)$ &$( 0_0, 5_0, 7_2)$
&$( 0_0, 8_2, 13_2)$ &$( 0_0, 10_0, 14_3)$
\end{tabular}}

{\footnotesize \noindent
\begin{tabular}{lllllllll}
$Q'$ &$( 0_0, 8_0, 3_2, 0_1,14_2)$ &$( 1_1, 6_2, 0_3,11_3,13_3)$
&$( 2_2, 7_3, 7_0, 7_1, 0_2)$ &$( 4_0, 8_1, 6_0, 9_1,14_3)$
\end{tabular}}

{\footnotesize \noindent
\begin{tabular}{lllllllllll}
$P$ &$( 3_1,13_3, 9_2)$ &$( 8_3, 7_1, 5_2)$ &$( 9_1, 5_3, 3_0)$
&$(11_2, 6_0, 7_0)$
&$( 2_2, 1_0, 5_0)$\\
&$(13_0,11_1,12_1)$ &$( 0_0, 3_3, 6_2)$ &$( 5_1, 2_1,10_0)$ &$(
1_1, 3_2,14_1)$
&$( 0_3, 8_2,14_0)$\\
&$( 8_0, 4_3, 1_3)$ &$(12_0, 9_0, 6_3)$ &$( 4_0, 2_0, 9_3)$ &$(
7_2, 0_2, 0_1)$
&$(11_3,12_2,10_3)$\\
&$( 7_3, 1_2, 4_1)$ &$(10_1, 4_2,14_3)$ &$(10_2,14_2,13_2)$
&$(13_1, 6_1,12_3)$ &$( 2_3, 8_1,11_0)$
\end{tabular}}

{\footnotesize \vspace{5pt}
\noindent {$(\alpha,\beta)=(24,5) : $}\\[5pt]
\noindent
\begin{tabular}{lllllllllll}
$P$
&$( 3_2,10_1,14_3)$
&$( 2_1, 6_0, 0_1)$
&$( 8_0,13_1, 3_1)$
&$( 6_2,12_3,11_1)$
&$( 0_0,14_0,10_3)$\\
&$( 1_1,12_2, 5_2)$
&$( 4_0, 8_1,12_1)$
&$( 2_2, 4_2, 9_3)$
&$(10_2, 0_2, 3_0)$
&$( 1_0, 8_3,11_2)$\\
&$( 6_1, 7_2, 7_1)$
&$( 9_1, 5_3, 6_3)$
&$(13_0, 1_2,11_0)$
&$(14_2, 1_3,10_0)$
&$( 5_1, 4_3, 5_0)$\\
&$( 3_3, 8_2, 9_2)$
&$( 0_3, 2_0, 7_0)$
&$(11_3, 9_0,13_3)$
&$(12_0, 4_1,13_2)$
&$( 7_3,14_1, 2_3)$\\
$Q_1$
&$( 0_0, 1_1, 8_2)$
&$( 3_3,13_0,11_1)$
&$( 6_2, 9_1, 4_3)$
&$( 8_0, 8_3, 4_2)$\\
$Q_2$
&$( 0_0, 3_3,12_1)$
&$( 1_1,10_0,14_0)$
&$( 2_2,11_1,13_2)$
&$( 6_2, 7_3,14_3)$\\
$Q_3$
&$( 0_0, 6_2, 7_0)$
&$( 1_1, 6_3,10_3)$
&$( 2_2, 8_0, 2_3)$
&$( 5_1,10_2,12_1)$\\
\end{tabular}}

{\footnotesize \noindent
\begin{tabular}{lllllllll}
$T'$
&$( 0_0, 2_2, 1_3, 13_2, 14_1)$
&$( 0_0, 12_3, 9_2, 11_0, 13_2)$
&$( 0_0, 1_3, 4_0, 2_2, 3_1)$
\end{tabular}}

{\footnotesize \vspace{5pt}
\noindent {$(\alpha,\beta)=(26,3) : $}\\[5pt]
\noindent
\begin{tabular}{lllllllllll}
$P$
&$( 3_3,10_2,11_0)$
&$( 5_1, 0_2, 9_3)$
&$( 2_2, 7_3,12_0)$
&$( 0_0, 9_0, 7_0)$
&$( 0_3,14_0,14_3)$\\
&$( 9_1, 2_0,11_1)$
&$( 6_2, 8_3,10_3)$
&$( 1_1,10_1, 6_0)$
&$( 4_0,12_3,13_2)$
&$( 8_0, 2_3, 6_3)$\\
&$(13_1, 7_2,14_1)$
&$( 6_1,13_0, 8_2)$
&$( 1_0, 4_2,12_2)$
&$(14_2, 7_1, 5_2)$
&$( 5_0,11_2, 9_2)$\\
&$( 2_1, 5_3, 3_0)$
&$( 3_2, 1_3, 4_1)$
&$(11_3, 1_2,12_1)$
&$( 4_3, 3_1,13_3)$
&$(10_0, 0_1, 8_1)$\\
$Q_1$
&$( 0_0, 4_0, 6_1)$
&$( 1_1, 6_2, 7_2)$
&$( 2_2, 9_3, 2_3)$
&$( 5_1, 5_0, 1_3)$\\
$Q_2$ &$( 0_0, 1_0, 6_3)$ &$( 1_1, 7_3, 1_2)$ &$( 2_2, 6_2, 5_0)$
&$( 5_1, 9_1,14_3)$\\
$S$ &$( 0_0, 1_1, 2_2)$ &$( 0_0, 4_3, 11_1)$ &$( 0_0, 5_0, 7_2)$
&$( 0_0, 8_2, 13_2)$ &$( 0_0, 10_0, 14_3)$
\end{tabular}}

{\footnotesize \noindent
\begin{tabular}{lllllllllll}
 $T'$ &$( 0_0, 3_3, 6_3, 9_3, 12_3)$ &$( 0_0,
12_0, 9_1,6_2, 3_1)$ &$( 0_0, 3_0,6_3, 9_0,12_1)$
\end{tabular}}

{\footnotesize \vspace{5pt}
\noindent {$(\alpha,\beta)=(27,2) : $}\\[5pt]
\noindent
\begin{tabular}{lllllllllll}
$P$
&$( 5_1, 8_0,12_0)$
&$( 0_2, 1_3,14_3)$
&$( 8_3,13_0, 6_0)$
&$(10_1, 9_2,11_0)$
&$( 4_3, 6_1,12_3)$\\
&$( 2_0, 1_2, 7_0)$
&$( 3_3, 4_0,13_3)$
&$( 9_1, 9_0, 4_2)$
&$( 0_0, 1_1, 2_2)$
&$(14_1, 9_3, 3_0)$\\
&$(13_1, 7_2,11_1)$
&$( 7_3, 3_1, 5_2)$
&$(11_3, 2_1,12_1)$
&$( 7_1, 0_1,13_2)$
&$( 0_3, 3_2,10_0)$\\
&$( 5_0,11_2, 2_3)$ &$( 6_2, 1_0,14_0)$ &$(10_2,14_2,10_3)$ &$(
5_3,12_2, 6_3)$ &$( 8_2, 4_1, 8_1)$
\end{tabular}}

{\footnotesize  \noindent
\begin{tabular}{lllllllllll}
$Q_1$ &$( 0_0, 3_3, 7_3)$ &$( 1_1, 2_1, 7_0)$ &$( 2_2, 5_0, 1_2)$
&$( 6_2, 9_1,11_3)$\\
$Q_2$ &$( 0_0, 5_1,14_0)$ &$( 1_1,13_0, 2_3)$ &$( 2_2, 9_3,12_1)$
&$( 6_2, 1_3,13_2)$\\
$Q_3$ &$( 0_0,11_3, 4_1)$ &$( 2_2, 1_0, 4_2)$
&$( 3_3, 2_0, 0_1)$
&$( 5_1,12_2,10_3)$\\
$Q_4$
&$( 0_0, 0_3, 7_1)$
&$( 2_2, 7_2,14_1)$
&$( 4_0, 6_1, 2_3)$
&$( 6_2, 2_0,10_3)$
\end{tabular}}

{\footnotesize \vspace{5pt}
\noindent {$(\alpha,\beta)=(28,1) : $}\\[5pt]
\noindent
\begin{tabular}{lllllllllll}
$P$

&$( 0_3,12_3,14_1)$
&$( 3_1,14_0,13_2)$
&$( 0_0, 1_1, 8_1)$
&$( 2_0, 5_2, 6_3)$
&$( 2_1, 7_2,11_0)$\\
&$( 0_2,12_2, 4_1)$
&$(11_3, 6_1,13_3)$
&$(12_0, 1_0,11_1)$
&$( 2_2,10_1, 9_2)$
&$( 5_1, 8_0, 9_3)$\\
&$( 9_0, 1_3, 7_1)$
&$(10_2,14_2,14_3)$
&$( 5_0,11_2, 3_0)$
&$( 4_2, 6_0, 0_1)$
&$(13_1, 4_3, 8_3)$\\
&$( 7_3,13_0,10_0)$
&$( 6_2, 5_3, 8_2)$
&$( 9_1, 3_2,12_1)$
&$( 3_3, 4_0, 2_3)$
&$( 1_2, 7_0,10_3)$\\
\end{tabular}}

{\footnotesize  \noindent
\begin{tabular}{lllllllllll}
$Q_1$ &$( 0_0, 2_2, 5_1)$ &$( 1_1,13_0,14_0)$ &$( 3_3, 7_2,10_3)$
&$( 6_2, 6_1,14_3)$\\
$Q_2$
&$( 0_0, 8_0,10_3)$
&$( 1_1, 1_0,12_3)$
&$( 2_2, 3_2, 8_3)$
&$( 5_1, 9_1, 7_2)$\\
$Q_3$
&$( 0_0,10_2, 7_1)$
&$( 2_2, 0_3,10_0)$
&$( 5_1, 6_1, 8_3)$
&$( 6_2, 5_0,13_3)$\\
$Q_4$
&$( 0_0, 0_3, 5_2)$
&$( 1_1,14_1,13_3)$
&$( 4_0, 0_1, 1_2)$
&$( 6_2, 8_3, 2_0)$
\end{tabular}}

\section*{Appendix D for Lemma~\ref{L37-42}}

{\footnotesize
\noindent {$(\alpha,\beta)=(1,19) : $}\\[5pt]
\noindent
\begin{tabular}{lllllllllll}
$P_1'$
&$( 0_0, 1_1, 2_2, 3_0, 4_1, 5_2, 6_0)$
&$( 7_1, 9_0,11_2, 8_2,10_1,12_0, 0_2)$
&$(13_1, 1_0, 5_1, 2_1, 6_2, 3_2, 7_0)$\\
&$( 4_0, 8_1,13_0,11_1, 0_1, 9_2, 1_2)$
&$(10_0, 2_0, 9_1,12_2, 8_0, 3_1,11_0)$
&$( 4_2,10_2, 5_0,12_1, 7_2, 6_1,13_2)$\\
$P_2'$
&$( 0_0, 3_0,10_1, 1_1, 7_1,11_2, 2_1)$
&$( 2_2,12_0, 4_1, 0_2, 5_2, 8_2, 1_0)$
&$( 6_0, 9_0, 3_2,13_1, 4_0,11_1, 5_1)$\\
&$( 6_2, 1_2,11_0,12_2, 6_1, 7_0, 4_2)$ &$( 8_1, 2_0,13_2,
3_1,12_1, 9_2, 7_2)$ &$(10_0, 8_0, 0_1,10_2,13_0, 5_0, 9_1)$
\end{tabular}}

{\footnotesize \noindent
\begin{tabular}{lllllllllll}
$P_3'$ &$( 0_0, 9_0, 4_0, 1_1,12_0,13_0, 1_0)$ &$(
2_2, 3_2, 6_1, 7_1, 5_0,13_1, 1_2)$
&$( 3_0,12_2, 4_1, 5_1, 3_1, 8_2, 0_1)$\\
&$( 5_2, 8_1,13_2, 7_0,12_1,10_0, 7_2)$
&$( 6_0,11_1,10_2, 2_1, 8_0,10_1, 4_2)$
&$(11_2, 9_2,11_0, 0_2, 9_1, 6_2, 2_0)$\\
\end{tabular}}

{\footnotesize \vspace{5pt}
\noindent {$(\alpha,\beta)=(2,18) : $}\\
\noindent
\begin{tabular}{lllllllllll}
$Q_1$ &$( 0_0, 3_0, 7_1)$ &$( 2_2, 5_2,10_1)$& $Q_2$ &$( 0_0, 4_1,
8_2)$ &$( 1_1, 5_2, 9_0)$
\end{tabular}}

{\footnotesize \noindent
\begin{tabular}{lllllllll}
$P_1'$
&$( 0_0, 5_2, 8_2, 1_1, 4_1, 7_1,12_0)$
&$( 2_2, 6_0, 3_0,10_1, 1_0, 6_2,11_2)$
&$( 9_0, 3_2,10_0, 0_2, 5_1, 0_1, 7_2)$\\
&$(13_1, 9_2, 3_1, 2_1, 6_1, 7_0, 5_0)$
&$( 4_0, 4_2,11_0, 1_2,13_2,12_2,10_2)$
&$( 8_1, 2_0,12_1,11_1, 9_1,13_0, 8_0)$\\
$P_2'$
&$( 0_0, 9_0, 9_2, 1_1, 0_2, 4_1,13_1)$
&$( 2_2, 1_0, 2_0, 3_2, 6_1,11_2, 1_2)$
&$( 3_0, 5_1, 9_1,12_0, 3_1, 5_2,13_0)$\\
&$( 6_0,10_0,13_2, 4_0, 5_0, 2_1, 0_1)$
&$( 7_1,12_2,12_1, 7_0,10_2, 8_1, 7_2)$
&$( 8_2,11_1,11_0, 6_2, 8_0,10_1, 4_2)$\\
$T'$
&$( 0_0, 1_1, 2_2, 3_0, 4_1, 5_2,6_0)$
&$( 0_0, 2_2, 4_1, 10_1, 12_0,1_1, 13_2)$
&$( 0_0, 11_2, 5_2,2_1, 13_0, 10_2,8_0)$\\
&$( 0_0,3_1,9_1, 1_1, 4_2, 6_1, 12_1)$
\end{tabular}}

{\footnotesize \vspace{5pt}
\noindent {$(\alpha,\beta)=(3,17) : $}\\[5pt]
\noindent
\begin{tabular}{lllllllllll}
$Q_1$ &$( 0_0, 1_1, 2_2)$ &$( 3_0, 4_1, 5_2)$& $Q_2$ &$( 0_0, 3_0,
7_1)$
&$( 2_2, 4_1,11_2)$\\
$Q_3$
&$( 0_0, 4_1, 8_2)$
&$( 1_1, 3_0,11_2)$
\end{tabular}}

{\footnotesize \noindent
\begin{tabular}{lllllllll}
$P_1'$
&$( 0_0, 5_2, 1_1, 4_1, 6_0, 2_2, 9_0)$
&$( 3_0, 8_2,11_2,12_0, 7_1, 0_2, 5_1)$
&$(10_1,13_1, 7_0, 2_1,10_0, 1_0,11_1)$\\
&$( 3_2, 6_2, 0_1, 4_0, 4_2, 5_0, 8_0)$
&$( 8_1, 3_1,12_1,12_2, 9_1, 7_2,11_0)$
&$( 9_2, 2_0,13_2,13_0, 6_1, 1_2,10_2)$\\
$P_2'$
&$( 0_0,13_1, 5_2, 1_0, 1_1,12_0,13_0)$
&$( 2_2, 3_2, 4_1, 0_2, 5_0, 8_2, 0_1)$
&$( 3_0,11_1,10_2, 4_0, 7_2, 5_1, 6_1)$\\
&$( 6_0, 3_1,11_2,10_0,12_1, 7_0, 4_2)$
&$( 7_1,12_2,11_0, 9_2, 9_1, 6_2, 2_0)$
&$( 9_0, 8_1, 8_0, 2_1,13_2,10_1, 1_2)$
\end{tabular}}

{\footnotesize \vspace{5pt}
\noindent {$(\alpha,\beta)=(4,16) : $}\\[5pt]
\noindent
\begin{tabular}{lllllllllll}
$Q_1$ &$( 0_0, 1_1, 2_2)$ &$( 3_0, 4_1, 5_2)$ &$Q_2$ &$( 0_0, 3_0,
7_1)$
&$( 2_2, 4_1,11_2)$\\
$Q_3$ &$( 0_0, 4_1, 8_2)$ &$( 1_1, 3_0,11_2)$ &$Q_4$ &$( 0_0, 5_2,
9_0)$ &$( 1_1, 4_1, 0_2)$
\end{tabular}}

{\footnotesize \noindent
\begin{tabular}{lllllllll}
$P_1'$
&$( 0_0, 6_0, 1_1, 5_2, 2_2, 3_0,10_1)$
&$( 4_1, 7_1,13_1, 8_2,12_0, 9_0, 1_0)$
&$(11_2, 0_2, 6_2,13_0, 3_2, 4_0, 2_1)$\\
&$( 5_1,12_2, 6_1, 9_2, 4_2, 7_0, 2_0)$
&$( 8_1, 7_2,13_2,11_1,12_1, 3_1,11_0)$
&$(10_0, 9_1, 1_2, 8_0, 0_1, 5_0,10_2)$\\
$P_2'$ &$( 0_0,11_2,13_1, 2_2,12_0, 0_1, 1_0)$ &$( 1_1, 4_0,
4_1,10_0, 5_2, 5_1, 5_0)$
&$( 3_0, 8_1, 3_1, 9_0,11_1,10_1, 1_2)$\\
&$( 6_0, 9_2,10_2, 7_0, 8_0, 6_2, 7_2)$
&$( 7_1, 2_0, 8_2,13_0,13_2, 2_1, 4_2)$
&$( 0_2, 6_1,12_1,12_2,11_0, 3_2, 9_1)$
\end{tabular}}

{\footnotesize \vspace{5pt}
\noindent {$(\alpha,\beta)=(5,15) : $}\\[5pt]
\noindent
\begin{tabular}{lllllllllll}
$P_1'$
&$( 0_0, 1_1, 2_2, 3_0, 4_1, 5_2, 8_2)$
&$( 6_0, 7_1,10_1,13_1,12_0, 9_0, 0_2)$
&$(11_2, 2_1, 9_2, 1_0, 4_0, 7_0,12_2)$\\
&$( 3_2, 8_1,13_0, 5_1,10_0, 5_0,12_1)$
&$( 6_2, 3_1,11_0, 6_1, 9_1, 0_1, 7_2)$
&$(11_1, 4_2,13_2, 2_0,10_2, 1_2, 8_0)$\\
$P_2'$
&$( 0_0,11_2, 2_2, 1_0, 6_0,13_1, 5_0)$
&$( 1_1,12_0,11_1, 8_2, 9_2, 3_0, 0_1)$
&$( 4_1,10_0, 7_2, 8_1,13_2, 4_0,12_2)$\\
&$( 5_2, 6_2, 9_1, 0_2, 8_0, 5_1, 2_0)$
&$( 7_1,13_0,12_1, 7_0, 6_1, 3_2, 4_2)$
&$( 9_0, 1_2,10_1,11_0, 2_1, 3_1,10_2)$
\end{tabular}}

{\footnotesize \vspace{5pt}
\noindent {$(\alpha,\beta)=(6,14) : $}\\[5pt]
\noindent
\begin{tabular}{lllllllllll}
$P_1$
&$( 0_0, 1_1, 2_2)$
&$( 3_0, 4_1, 5_2)$
&$( 6_0, 7_1,10_1)$
&$( 8_2,11_2,13_1)$
&$( 9_0,12_0, 0_2)$\\
&$( 1_0, 3_2, 7_0)$
&$( 2_1, 4_0, 8_1)$
&$( 5_1, 6_2, 9_2)$
&$(10_0,13_0, 1_2)$
&$(11_1, 0_1, 4_2)$\\
&$(12_2, 7_2,13_2)$
&$( 2_0, 6_1,11_0)$
&$( 3_1, 9_1,10_2)$
&$( 5_0, 8_0,12_1)$
\end{tabular}}

{\footnotesize \noindent
\begin{tabular}{lllllllll}
$P_1'$
&$( 0_0, 5_2,10_1, 1_1, 8_2, 2_1, 7_1)$
&$( 2_2, 9_0, 4_1,12_0, 3_0,11_2, 3_2)$
&$( 6_0,13_1, 5_1,11_1, 0_2,10_0, 6_1)$\\
&$( 1_0, 0_1, 8_0, 7_0, 7_2, 4_0, 3_1)$ &$( 6_2,
1_2,12_1,12_2,11_0, 9_2, 5_0)$ &$( 8_1, 4_2,13_2, 2_0,10_2,13_0,
9_1)$\\

 $P_2'$ &$( 0_0,13_1, 5_2, 1_0, 1_1, 4_0, 4_2)$ &$( 2_2, 2_1,
1_2, 9_0, 0_1, 8_2, 5_0)$
&$( 3_0,11_1, 8_0, 3_2, 2_0,10_1,12_2)$\\
&$( 4_1, 5_1,11_0,10_0,10_2, 7_1, 9_2)$
&$( 6_0, 8_1,13_2, 7_0,12_1, 0_2, 3_1)$
&$(11_2,13_0, 7_2, 6_2, 9_1,12_0, 6_1)$
\end{tabular}}

{\footnotesize \vspace{5pt}
\noindent {$(\alpha,\beta)=(7,13) : $}\\[5pt]
\noindent
\begin{tabular}{lllllllllll}
$P_1$
&$( 0_0, 1_1, 2_2)$
&$( 3_0, 4_1, 5_2)$
&$( 6_0, 7_1,10_1)$
&$( 8_2,11_2,13_1)$
&$( 9_0,12_0, 0_2)$\\
&$( 1_0, 3_2, 7_0)$
&$( 2_1, 4_0, 8_1)$
&$( 5_1, 6_2, 9_2)$
&$(10_0,13_0, 1_2)$
&$(11_1, 0_1, 4_2)$\\
&$(12_2, 7_2,13_2)$
&$( 2_0, 6_1,11_0)$
&$( 3_1, 9_1,10_2)$
&$( 5_0, 8_0,12_1)$
\end{tabular}}

{\footnotesize \noindent
\begin{tabular}{lllllllll}
$P_1'$
&$( 0_0, 5_2,10_1, 1_1, 8_2, 2_1, 7_1)$
&$( 2_2, 9_0, 4_1,12_0, 3_0,11_2, 3_2)$
&$( 6_0,13_1, 5_1,11_1, 0_2, 8_1, 7_2)$\\
&$( 1_0,13_0, 8_0, 7_0, 3_1,13_2, 4_2)$ &$( 4_0, 0_1,12_1,
1_2,11_0,10_0, 9_1)$ &$( 6_2, 2_0,10_2,12_2, 6_1, 9_2, 5_0)$
\end{tabular}}

{\footnotesize \noindent
\begin{tabular}{lllllllll}
$P_2'$ &$( 0_0, 2_1, 1_1, 3_2, 4_1, 1_0, 6_1)$ &$( 2_2, 0_2, 5_0,
7_1,13_0,10_2, 7_0)$
&$( 3_0, 6_2, 9_1,12_0, 3_1, 5_2, 8_1)$\\
&$( 6_0,12_2,11_0, 5_1, 2_0,10_1, 4_2)$
&$( 8_2, 0_1,13_1, 1_2, 9_0,11_1,12_1)$
&$(11_2, 9_2, 8_0,10_0, 7_2, 4_0,13_2)$
\end{tabular}}

{\footnotesize \vspace{5pt}
\noindent {$(\alpha,\beta)=(8,12) : $}\\[5pt]
\noindent
\begin{tabular}{lllllllllll}
$T$
&$( 0_0, 1_1, 5_2)$\\
\end{tabular}}

{\footnotesize \noindent
\begin{tabular}{lllllllll}
$P_1'$
&$( 0_0, 2_2, 4_1, 1_1, 3_0, 5_2, 7_1)$
&$( 6_0, 8_2, 0_2,11_2,13_1,10_1, 2_1)$
&$( 9_0,12_0, 1_0, 8_1, 3_1, 4_0,10_0)$\\
&$( 3_2,12_2, 5_1,11_1, 4_2, 6_2, 5_0)$
&$( 7_0,13_0, 8_0, 2_0, 9_1, 6_1,12_1)$
&$( 9_2, 7_2,10_2, 0_1,13_2, 1_2,11_0)$\\
$P_2'$
&$( 0_0, 9_0, 4_0, 3_0,10_1, 6_2,13_0)$
&$( 1_1, 0_2, 4_1, 5_1, 0_1,12_0,11_1)$
&$( 2_2, 7_0, 5_0, 6_0,11_0,13_1, 4_2)$\\
&$( 5_2,10_0, 9_1, 8_1,13_2, 3_2, 1_2)$
&$( 7_1, 9_2, 8_0,12_2,10_2, 1_0, 3_1)$
&$( 8_2, 6_1, 2_1,12_1, 2_0,11_2, 7_2)$
\end{tabular}}

{\footnotesize \vspace{5pt}
\noindent {$(\alpha,\beta)=(9,11) : $}\\
\noindent
\begin{tabular}{lllllllllll}
$T$
&$( 0_0, 2_2, 13_1)$
&$( 0_0, 7_1, 3_2)$
&$( 0_0, 8_2, 8_1)$
\end{tabular}}

{\footnotesize \noindent
\begin{tabular}{lllllllllll}
$P_1'$
&$( 0_0, 1_1, 2_2, 3_0, 4_1, 5_2, 6_0)$
&$( 7_1,10_1,13_1, 8_2,11_2, 0_2, 9_2)$
&$( 9_0,12_0, 1_0, 5_1,10_0, 6_2,11_1)$\\
&$( 2_1, 3_1, 9_1, 4_2, 0_1, 5_0,11_0)$
&$( 3_2, 7_0, 8_0,13_0, 2_0, 6_1,12_1)$
&$( 4_0, 8_1,10_2, 1_2, 7_2,12_2,13_2)$\\
$Q'$
&$( 0_0, 9_2, 0_1, 5_2,10_2,11_2,13_0)$
&$( 1_1, 6_2, 8_0, 9_0, 4_0, 5_0, 9_1)$
&$( 3_0, 8_2, 3_1,11_1, 6_1, 7_2, 5_1)$
\end{tabular}}

{\footnotesize \vspace{5pt}
\noindent {$(\alpha,\beta)=(10,10) : $}\\[5pt]
\noindent
\begin{tabular}{lllllllllll}
$Q_1$ &$( 0_0, 1_1, 2_2)$ &$( 3_0, 4_1, 5_2)$ &$Q_2$ &$( 0_0, 3_0,
7_1)$
&$( 2_2, 4_1,11_2)$\\
$Q_3$ &$( 0_0, 4_1, 9_0)$ &$( 1_1, 5_2, 8_2)$ &$Q_4$ &$( 0_0,
5_2,10_1)$
&$( 1_1, 3_0,10_2)$\\
$Q_5$
&$( 0_0,13_1, 2_1)$
&$( 2_2, 3_0, 1_2)$
\end{tabular}}

{\footnotesize \noindent
\begin{tabular}{lllllllll}
$P_1'$
&$( 0_0, 1_0, 2_2, 5_2, 6_0, 1_1,13_0)$
&$( 3_0, 8_2, 4_1, 7_1, 6_2,11_0, 7_2)$
&$( 9_0, 4_0, 9_1, 0_2,10_0, 3_2, 5_1)$\\
&$(10_1, 1_2,13_1, 9_2, 4_2, 8_0, 5_0)$
&$(11_2, 2_0, 0_1, 7_0, 6_1,11_1,12_1)$
&$(12_0,10_2, 8_1, 3_1, 2_1,12_2,13_2)$
\end{tabular}}

{\footnotesize \vspace{5pt}
\noindent {$(\alpha,\beta)=(11,9) : $}\\[5pt]
\noindent
\begin{tabular}{lllllllllll}
$P_1$
&$( 0_0, 1_1, 2_2)$
&$( 3_0, 4_1, 5_2)$
&$( 6_0, 7_1,10_1)$
&$( 8_2,11_2,13_1)$
&$( 9_0,12_0, 0_2)$\\
&$( 1_0, 3_2, 8_1)$
&$( 2_1, 4_0, 9_2)$
&$( 5_1, 6_2, 1_2)$
&$( 7_0,12_2, 3_1)$
&$(10_0, 0_1, 9_1)$\\
&$(11_1,10_2,13_2)$
&$(13_0, 2_0, 6_1)$
&$( 4_2, 7_2,11_0)$
&$( 5_0, 8_0,12_1)$
\end{tabular}}

{\footnotesize \noindent
\begin{tabular}{lllllllll}
$P_1'$
&$( 0_0, 7_1, 6_2, 1_1,10_1, 3_0, 5_1)$
&$( 2_2,11_2, 6_0, 8_2, 7_2,11_1, 6_1)$
&$( 4_1, 2_0,12_2,11_0, 1_2, 2_1, 3_1)$\\
&$( 5_2,10_0, 9_2, 5_0, 4_0, 8_1,10_2)$
&$( 9_0, 4_2,13_1,12_1,13_0, 3_2, 8_0)$
&$(12_0, 7_0, 9_1, 0_2,13_2, 1_0, 0_1)$
\end{tabular}}

{\footnotesize \vspace{5pt}
\noindent {$(\alpha,\beta)=(12,8) : $}\\[5pt]
\noindent
\begin{tabular}{lllllllllll}
$P_1$
&$( 0_0, 1_1, 2_2)$
&$( 3_0, 4_1, 5_2)$
&$( 6_0, 7_1,10_1)$
&$( 8_2,11_2,13_1)$
&$( 9_0,12_0, 0_2)$\\
&$( 1_0, 3_2, 7_0)$
&$( 2_1, 4_0, 8_1)$
&$( 5_1, 6_2, 9_2)$
&$(10_0,13_0, 1_2)$
&$(11_1, 0_1, 4_2)$\\
&$(12_2, 7_2,13_2)$ &$( 2_0, 6_1,11_0)$ &$( 3_1, 9_1,10_2)$ &$(
5_0, 8_0,12_1)$\\

$P_2$ &$( 0_0, 5_2,10_1)$ &$( 1_1, 8_2, 2_1)$ &$( 2_2, 9_0, 3_2)$
&$( 3_0,11_2, 5_1)$
&$( 4_1,13_1, 6_2)$\\
&$( 6_0, 1_0, 1_2)$
&$( 7_1,13_0, 7_2)$
&$(12_0, 9_2, 9_1)$
&$( 0_2,11_1, 2_0)$
&$( 4_0, 3_1,13_2)$\\
&$( 7_0, 0_1, 8_0)$
&$( 8_1, 4_2,10_2)$
&$(10_0, 5_0,11_0)$
&$(12_2, 6_1,12_1)$
\end{tabular}}

{\footnotesize \noindent
\begin{tabular}{lllllllll}
$P_1'$
&$( 0_0, 9_0, 5_1, 4_1, 0_2, 1_1, 9_2)$
&$( 2_2,13_1,13_0, 5_2, 7_0,12_1, 1_0)$
&$( 3_0, 2_1, 4_2, 7_1, 6_1, 8_2, 0_1)$\\
&$( 6_0, 3_2, 3_1,12_0,13_2,10_0, 7_2)$
&$(10_1, 1_2,10_2, 4_0, 9_1, 6_2, 2_0)$
&$(11_2,11_1, 8_0,12_2,11_0, 8_1, 5_0)$
\end{tabular}}

{\footnotesize \vspace{5pt}
\noindent {$(\alpha,\beta)=(13,7) : $}\\[5pt]
\noindent
\begin{tabular}{lllllllllll}
$P_1$
&$( 0_0, 1_1, 2_2)$
&$( 3_0, 4_1, 5_2)$
&$( 6_0, 7_1,10_1)$
&$( 8_2,11_2,13_1)$
&$( 9_0,12_0, 0_2)$\\
&$( 1_0, 3_2, 7_0)$
&$( 2_1, 4_0, 8_1)$
&$( 5_1, 6_2, 9_2)$
&$(10_0,13_0, 1_2)$
&$(11_1, 0_1, 4_2)$\\
&$(12_2, 7_2,13_2)$ &$( 2_0, 6_1,11_0)$ &$( 3_1, 9_1,10_2)$ &$(
5_0, 8_0,12_1)$
\end{tabular}}

{\footnotesize \noindent
\begin{tabular}{lllllllll}
$P_2$ &$( 0_0, 5_2,10_1)$ &$( 1_1, 8_2, 2_1)$ &$( 2_2, 9_0, 3_2)$
&$( 3_0,11_2, 5_1)$
&$( 4_1,13_1, 6_2)$\\
&$( 6_0,12_0, 1_2)$
&$( 7_1, 9_2, 6_1)$
&$( 0_2,12_2, 9_1)$
&$( 1_0,13_0, 7_2)$
&$( 4_0, 3_1,13_2)$\\
&$( 7_0, 0_1, 8_0)$
&$( 8_1, 4_2,10_2)$
&$(10_0, 5_0,11_0)$
&$(11_1, 2_0,12_1)$
\end{tabular}}

{\footnotesize \noindent
\begin{tabular}{lllllllll}
$P_1'$
&$( 0_0, 9_0, 7_0, 4_1, 0_2, 1_1, 3_1)$
&$( 2_2, 1_0, 5_2, 4_0, 1_2,10_2, 5_1)$
&$( 3_0, 9_2, 6_0, 7_2,13_1, 8_1, 5_0)$\\
&$( 7_1,10_0, 9_1,11_1, 8_2, 6_1,12_2)$
&$(10_1,13_0,12_1, 6_2, 2_0,11_2,13_2)$
&$(12_0, 0_1,11_0, 3_2, 8_0, 2_1, 4_2)$
\end{tabular}}

{\footnotesize \vspace{5pt}
\noindent {$(\alpha,\beta)=(15,5) : $}\\[5pt]
\noindent
\begin{tabular}{lllllllllll}
$P_1$
&$( 0_0, 1_1, 4_1)$
&$( 2_2, 3_0, 6_0)$
&$( 5_2, 8_2, 9_0)$
&$( 7_1,10_1,12_2)$
&$(11_2, 1_0, 6_2)$\\
&$(12_0, 4_0, 7_0)$
&$(13_1, 0_2, 5_1)$
&$( 2_1, 3_1, 7_2)$
&$( 3_2, 9_2, 4_2)$
&$( 8_1, 0_1, 9_1)$\\
&$(10_0, 1_2, 5_0)$
&$(11_1,12_1,13_2)$
&$(13_0, 2_0, 8_0)$
&$( 6_1,10_2,11_0)$\\
$T$
&$( 0_0, 2_2, 13_1)$
&$( 0_0, 7_1, 3_2)$
&$( 0_0, 8_2, 8_1)$
\end{tabular}}

{\footnotesize \noindent
\begin{tabular}{lllllllll}
$Q'$
&$( 0_0, 6_0,11_1, 9_0, 0_2,13_2,12_2)$
&$( 1_1,10_0, 5_1, 2_1,11_0, 2_2, 0_1)$
&$( 8_2, 3_2,11_2,12_0, 3_1, 1_0, 6_1)$
\end{tabular}}

{\footnotesize \vspace{5pt}
\noindent {$(\alpha,\beta)=(16,4) : $}\\[5pt]
\noindent
\begin{tabular}{lllllllll}
$Q'$ &$( 0_0, 9_0, 7_2, 1_0, 6_0, 0_1,12_0)$ &$( 1_1,10_1,
5_1,11_1,13_2, 4_0,9_1)$ &$( 2_2,10_0, 5_2,13_1,1_2, 3_2, 4_2)$
\end{tabular}

\noindent
\begin{tabular}{lllllllllll}
$P_1$
&$( 0_0, 3_0, 6_0)$
&$( 1_1, 4_1, 7_1)$
&$( 2_2, 5_2, 0_2)$
&$( 8_2, 6_2, 7_2)$
&$( 9_0, 0_1, 8_0)$\\
&$(10_1,12_2, 1_2)$
&$(11_2, 5_1,13_2)$
&$(12_0,10_0, 4_2)$
&$(13_1, 7_0,12_1)$
&$( 1_0, 6_1, 9_1)$\\
&$( 2_1,11_1, 3_1)$
&$( 3_2, 9_2,11_0)$
&$( 4_0,13_0, 5_0)$
&$( 8_1, 2_0,10_2)$\\
$T$
&$( 0_0, 2_2, 13_1)$
&$( 0_0, 7_1, 3_2)$
&$( 0_0, 1_1, 5_2)$
\end{tabular}}

{\footnotesize \vspace{5pt}
\noindent {$(\alpha,\beta)=(17,3) : $}\\[5pt]
\noindent
\begin{tabular}{lllllllllll}
$P_1$
&$( 0_0, 1_1, 2_2)$
&$( 3_0, 6_0,10_1)$
&$( 4_1, 7_1, 9_0)$
&$( 5_2, 8_2,12_0)$
&$(11_2,13_1, 0_1)$\\
&$( 0_2,11_1,13_2)$
&$( 1_0, 4_2, 9_1)$
&$( 2_1, 6_2, 5_0)$
&$( 3_2,13_0, 8_0)$
&$( 4_0,12_2, 6_1)$\\
&$( 5_1, 1_2,10_2)$
&$( 7_0, 3_1,12_1)$
&$( 8_1, 2_0, 7_2)$
&$( 9_2,10_0,11_0)$\\
$S$
 &$( 0_0,11_2, 3_1)$
 &$( 0_0,13_1, 6_2)$
\end{tabular}}

{\footnotesize \vspace{5pt}
\noindent {$(\alpha,\beta)=(18,2) : $}\\[5pt]
\noindent
\begin{tabular}{lllllllllll}
$P_1$
&$( 0_0, 1_1, 2_2)$
&$( 3_0, 6_0,10_1)$
&$( 4_1, 7_1, 9_0)$
&$( 5_2, 8_2,12_0)$
&$(11_2,13_1, 6_2)$\\
&$( 0_2, 1_2, 9_1)$
&$( 1_0, 7_0, 6_1)$
&$( 2_1, 5_0,13_2)$
&$( 3_2, 0_1, 8_0)$
&$( 4_0,13_0,10_2)$\\
&$( 5_1, 2_0, 7_2)$
&$( 8_1,12_2, 4_2)$
&$( 9_2,10_0,11_0)$
&$(11_1, 3_1,12_1)$\\
$S$
  &$( 0_0,10_1, 4_2)$
  &$( 0_0, 2_1, 1_2)$
  &$( 0_0, 3_2, 8_1)$
\end{tabular}}

{\footnotesize \vspace{5pt}
\noindent {$(\alpha,\beta)=(19,1) : $}\\[5pt]
\noindent
\begin{tabular}{lllllllllll}
$P_1$
&$( 0_0, 3_0, 7_1)$
&$( 1_1, 4_1, 8_2)$
&$( 2_2, 5_2, 9_0)$
&$( 6_0,11_2, 9_2)$
&$(10_1, 7_2,13_2)$\\
&$(12_0, 6_2,11_0)$
&$(13_1, 7_0, 2_0)$
&$( 0_2,10_0, 8_0)$
&$( 1_0,12_2, 6_1)$
&$( 2_1,11_1, 3_1)$\\
&$( 3_2, 4_2, 9_1)$
&$( 4_0, 0_1,12_1)$
&$( 5_1, 1_2,10_2)$
&$( 8_1,13_0, 5_0)$
\end{tabular}}

\section*{Appendix E for Lemma~\ref{RGDD84}}

{\footnotesize
\noindent {$(\alpha,\beta)=(0,31) : $ \ \  $I=\{(i_0,(i+17)_2),(i_1,(i+17)_3)| i \in Z_{21}\}$.}\\[5pt]
\noindent
\begin{tabular}{lllllllllll}
$Q_1'$
&$( 0_0, 2_2, 4_0, 1_1, 3_3, 5_1,10_2)$
&$( 6_2, 9_1,15_3, 8_0,11_3,14_2,19_3)$
&$( 7_3,12_0,17_1, 3_0,13_1, 1_2,16_0)$\\
&$(18_2, 4_1, 6_0, 2_3, 7_1, 5_2,13_3)$\\
$Q_2'$
&$( 0_0, 9_1, 3_3,10_2,16_0, 1_1,11_3)$
&$( 2_2,12_0, 4_1, 4_0,18_2, 5_1,19_3)$
&$( 6_2,15_3, 5_2, 7_3,17_1,15_0, 9_3)$\\
&$(13_1, 0_2,20_0, 7_1, 3_0, 8_2,20_3)$\\
\end{tabular}}

{\footnotesize  \noindent
\begin{tabular}{lllllllllll}
$Q_3'$ &$( 0_0,17_1, 2_2,13_1,10_3, 1_1,18_2)$ &$( 4_0, 2_3,15_0,
5_1, 5_2,11_3,13_2)$
&$( 7_3, 3_0,10_2,19_0,19_3,16_1, 7_2)$\\
&$(15_3, 2_0, 1_2,11_1,20_0,13_3,14_1)$\\
$Q_4'$
&$( 0_0, 2_3,11_1, 4_0,14_3, 9_1, 4_2)$
&$( 1_1, 5_2,10_0,13_1,13_3, 8_0, 5_3)$
&$( 2_2, 1_3,15_2,12_0, 7_1,10_2, 3_1)$\\
&$( 3_3, 2_0,11_3, 0_2,19_1,20_0,20_2)$\\
$T'$
&$( 0_0, 1_1, 2_2, 3_3, 4_0,19_1,13_2)$
&$( 0_0, 6_3,19_1,18_0, 3_3, 9_2,15_1)$
&$( 0_0, 8_2,16_0, 3_2,11_0,19_2,20_3)$
\end{tabular}}

{\footnotesize \vspace{5pt}
\noindent {$(\alpha,\beta)=(9,22) : $ \ \ $I=\{(i_0,(i+8)_2),(i_1,(i+4)_3)| i \in Z_{21}\}$.}\\[5pt]
\noindent
\begin{tabular}{lllllllllll}
$Q_1$ &$( 0_0, 6_2,13_1)$ &$( 2_2, 8_0,15_3)$ &$( 4_0,19_3, 3_1)$
&$( 5_1,11_3, 1_2)$\\
$Q_2$ &$( 0_0, 9_1,18_2)$ &$( 1_1, 8_0,19_3)$ &$( 2_2,11_3,17_1)$
&$( 3_3,10_2,19_0)$\\
$Q_3$ &$( 0_0,10_2, 0_1)$ &$( 1_1,11_3, 7_0)$ &$( 2_2,19_3, 8_1)$
&$( 3_3,18_2, 2_0)$
\end{tabular}}

{\footnotesize \noindent
\begin{tabular}{lllllllll}
$Q_1'$
&$( 0_0,14_2, 1_1,12_0, 2_2,16_0,12_1)$
&$( 3_3,17_1,11_0, 9_1, 2_3,15_0, 3_2)$
&$( 6_2, 3_0,15_3, 4_2, 7_3, 4_1, 4_3)$\\
&$(13_1, 5_3, 0_1,12_2,20_0,20_3, 8_2)$\\
$Q_2'$
&$( 0_0, 2_3, 8_2, 1_1,10_3, 4_0,14_3)$
&$( 2_2,12_1, 5_0, 0_1,10_0,10_2, 6_0)$
&$( 6_2, 1_3,10_1,18_3, 9_0,14_2,13_3)$\\
&$( 8_0, 5_3, 9_1,12_2,13_1,11_2,11_1)$\\

$T'$ &$( 0_0, 1_1,19_2, 6_3,10_0,11_2,16_3)$ &$( 0_0, 2_2, 4_3,
6_1, 3_2, 5_3, 1_2)$
&$( 0_0,17_3,18_0,15_1,16_3,20_0,19_2)$\\
&$( 0_0,13_3,15_0,11_3, 3_2, 5_0, 2_1)$
&$( 0_0, 4_1, 5_2,10_3,15_0, 2_1,20_2)$
&$( 0_0, 5_1, 4_3, 2_1, 3_2, 1_1,20_3)$\\
&$( 0_0, 3_3, 2_1, 1_0, 6_1,11_2,19_3)$
&$( 0_0, 8_1,10_3, 2_2,19_1,20_2,18_1)$
\end{tabular}}

{\footnotesize \vspace{5pt}
\noindent {$(\alpha,\beta)=(18,13) : $ \ \ $I=\{(i_0,(i+3)_1),(i_2,(i+12)_3)| i \in Z_{21}\}$.}\\[5pt]
\noindent
\begin{tabular}{lllllllllll}
$Q_1$
&$( 0_0, 2_2, 5_1)$
&$( 1_1, 3_3, 6_2)$
&$( 4_0, 7_3, 0_1)$
&$( 8_0, 1_2, 5_3)$\\
$Q_2$
&$( 0_0, 7_3, 9_1)$
&$( 1_1, 4_0, 1_2)$
&$( 2_2,11_3,20_0)$
&$( 3_3,18_2, 8_1)$\\
$Q_3$
&$( 0_0,10_2,15_3)$
&$( 1_1, 8_0,19_3)$
&$( 2_2, 4_0,12_1)$
&$( 5_1, 9_2, 5_3)$\\
$Q_4$ &$( 0_0,19_3, 5_2)$ &$( 1_1,10_2,20_0)$ &$( 3_3, 0_1,10_0)$
&$( 5_1,14_3,12_2)$\\
$Q_5$ &$( 0_0, 0_1,10_3)$ &$( 1_1,15_3,17_2)$ &$( 4_0, 8_1, 4_2)$
&$( 6_2,14_3, 5_0)$
\end{tabular}}

{\footnotesize  \noindent
\begin{tabular}{lllllllll}
$Q_6$ &$( 0_0, 2_3,12_1)$ &$( 1_1,16_0, 4_2)$ &$( 2_2,17_1,18_3)$
&$( 6_2,20_0,16_3)$
\end{tabular}}

{\footnotesize  \noindent
\begin{tabular}{lllllllll}
$Q_1'$
&$( 0_0,16_1, 7_3,12_0, 1_1, 3_0,20_2)$
&$( 2_2, 0_1,14_2, 2_0, 6_3,10_1, 9_3)$
&$( 3_3,20_1,11_2, 1_3, 1_0,17_2,19_1)$\\
&$( 4_0, 8_2,11_3, 6_0,19_3,19_2,11_1)$\\
$T'$
&$( 0_0,13_2,12_3,11_0, 3_2, 9_1, 8_2)$
&$( 0_0, 1_1, 9_0,17_2,11_3,19_2,20_3)$
&$( 0_0,20_1, 5_0,11_3, 3_0,16_1, 8_3)$\\
&$( 0_0, 6_3,19_0,18_3,17_2,16_1, 1_3)$
&$( 0_0, 6_2,12_0,18_2, 3_0, 9_2,15_1)$
&$( 0_0,13_1, 5_3,18_1,17_2,16_3,15_2)$
\end{tabular}}

{\footnotesize \vspace{5pt}
\noindent {$(\alpha,\beta)=(24,7) : $} \ \ $I=\{(i_0,(i+1)_3),(i_1,(i+20)_2)| i \in Z_{21}\}$.\\[5pt]
\noindent
\begin{tabular}{lllllllllll}
$Q_1$
&$( 0_0, 1_1, 2_2)$
&$( 3_3, 4_0, 9_1)$
&$( 5_1, 7_3,10_2)$
&$( 6_2, 8_0,11_3)$\\
$Q_2$
&$( 0_0, 6_2, 7_3)$
&$( 1_1, 4_0,14_2)$
&$( 3_3, 5_1,20_0)$
&$( 9_1, 2_3,13_2)$\\
$Q_3$
&$( 0_0, 9_1,15_3)$
&$( 1_1, 8_0, 1_2)$
&$( 2_2, 5_1, 2_3)$
&$( 4_0, 9_2, 1_3)$\\
$Q_4$
&$( 0_0,11_3,18_2)$
&$( 1_1,10_2,20_0)$
&$( 2_2, 9_1,19_3)$
&$( 3_3,16_0, 8_1)$\\
$Q_5$
&$( 0_0,17_1,13_2)$
&$( 1_1, 2_3, 7_0)$
&$( 2_2, 8_0, 6_3)$
&$( 6_2, 0_1,13_3)$\\
$Q_6$
&$( 0_0, 0_1, 5_3)$
&$( 1_1, 9_2,18_3)$
&$( 2_2,16_0,20_1)$
&$( 7_3, 1_2,14_0)$\\
$Q_7$
&$( 0_0, 1_2,12_1)$
&$( 1_1,10_3,19_0)$
&$( 2_2,20_0, 5_3)$
&$( 3_3,17_1,12_2)$\\
$Q_8$
&$( 0_0, 2_3,20_1)$
&$( 1_1, 1_3,20_2)$
&$( 3_3, 7_0,16_2)$
&$( 6_2,12_1, 2_0)$\\
\end{tabular}}

{\footnotesize  \noindent
\begin{tabular}{lllllllll}
$Q_1'$
&$( 0_0, 8_1, 7_3,17_1,19_0, 5_1,17_2)$
&$( 2_2,11_0,20_3, 1_2,20_1, 3_3,16_1)$\\
&$( 6_2, 5_3,11_2, 0_1,10_0, 5_2, 6_0)$
&$( 8_0,18_3, 7_2, 2_3,11_1,16_0, 8_3)$
\end{tabular}}

\section*{Appendix F for Lemma~\ref{L37-84}}

{\footnotesize
\noindent {$(\alpha,\beta)=(36,5) : $}\\[5pt]
\noindent
\begin{tabular}{lllllllllll}
$P$
    &$( 0_0,16_2, 6_3)$
    &$( 3_0,10_0, 2_2)$
    &$(19_1, 8_2,12_2)$
    &$( 4_0,16_1,10_2)$
    &$(12_1, 5_3,16_3)$
    &$( 9_2,17_2,12_3)$\\
    &$( 7_0,11_0,15_1)$
    &$(20_0,18_1, 0_3)$
    &$( 2_1,10_1,18_3)$
    &$( 0_2, 2_3,15_3)$
    &$( 5_1, 7_1,18_2)$
    &$( 5_0,15_0, 3_3)$\\
    &$(13_0,13_2,15_2)$
    &$(17_0, 7_3, 9_3)$
    &$(16_0,13_1, 4_2)$
    &$( 1_0, 3_1, 4_3)$
    &$(11_1, 1_3, 8_3)$
    &$(18_0, 1_2,11_2)$\\
    &$(12_0,19_2,19_3)$
    &$( 6_0,14_0,20_1)$
    &$( 2_0, 9_1, 5_2)$
    &$( 0_1, 4_1,14_1)$
    &$( 8_0,17_1,20_2)$
    &$( 7_2,14_2,13_3)$\\
    &$( 6_1, 3_2,11_3)$
    &$( 1_1, 6_2,20_3)$
    &$( 9_0, 8_1,17_3)$
    &$(19_0,10_3,14_3)$\\
    \end{tabular}}

{\footnotesize  \noindent
\begin{tabular}{lllllllll}
$Q_1$
    &$( 1_1, 7_2,16_3)$
    &$( 1_0, 2_1, 6_2)$
    &$( 0_1, 8_2, 0_3)$
    &$( 0_0, 2_0, 2_3)$\\
$Q_2$
    &$( 0_0,16_1,18_2)$
    &$( 1_0,16_2, 5_3)$
    &$( 2_0, 2_1,19_3)$
    &$( 0_1,20_2, 6_3)$\\
$Q_3$
    &$( 2_0,15_1, 1_3)$
    &$( 0_0,10_2,15_3)$
    &$( 1_0,16_1, 2_2)$
    &$( 2_1,18_2,14_3)$\\
$Q_4$
    &$( 2_0,19_2,16_3)$
    &$( 0_0, 3_1, 5_3)$
    &$( 1_0,11_1, 9_2)$
    &$( 1_1, 2_2, 0_3)$\\
$Q_5$
    &$( 0_0,11_1,11_2)$
    &$( 1_0, 6_1,19_3)$
    &$( 2_0, 0_2,12_3)$
    &$( 1_1,10_2,11_3)$\\
 \end{tabular}}

{\footnotesize \noindent
\begin{tabular}{lllllllllll}
$T'$
&$( 0_0, 1_1, 2_2, 3_3, 4_0, 5_1, 6_2)$
&$( 0_0, 3_3, 6_2, 1_1, 4_0,19_3,16_0)$
&$( 0_0, 5_1, 8_0, 3_3,18_2, 2_2,20_0)$\\
&$( 0_0,15_3, 9_1, 3_3,19_3,13_1,18_2)$
\end{tabular}}

{\footnotesize \vspace{5pt}
\noindent {$(\alpha,\beta)=(37,4) : $}\\[5pt]
\noindent
\begin{tabular}{lllllllllll}
$P$
    &$( 0_0,19_2, 0_3)$
    &$( 5_1,15_1, 9_3)$
    &$( 3_0,14_1,12_3)$
    &$( 0_1, 8_1,17_1)$
    &$( 6_0,14_0,19_1)$
    &$(13_2,10_3,20_3)$\\
    &$( 5_0, 9_1, 5_2)$
    &$(15_0, 8_2,12_2)$
    &$( 2_2,11_2,11_3)$
    &$( 3_2,16_2, 1_3)$
    &$( 7_0,17_0, 4_1)$
    &$(20_0, 6_1,13_3)$\\
    &$(16_1, 1_2,15_3)$
    &$(12_0,13_1,16_3)$
    &$(16_0, 6_3,19_3)$
    &$( 2_0,18_2, 7_3)$
    &$(18_1,14_3,18_3)$
    &$( 1_1,20_1,20_2)$\\
    &$( 3_1, 0_2, 8_3)$
    &$(11_1, 2_3, 4_3)$
    &$(13_0, 2_1,15_2)$
    &$( 4_0, 8_0, 3_3)$
    &$( 9_0,18_0,10_2)$
    &$( 7_1, 6_2,17_2)$\\
    &$( 1_0,10_1, 4_2)$
    &$(10_0,12_1,14_2)$
    &$(19_0, 7_2, 9_2)$
    &$(11_0, 5_3,17_3)$\\
    \end{tabular}}

{\footnotesize  \noindent
\begin{tabular}{lllllllll}
$Q_1$
    &$( 0_1, 7_2, 6_3)$
    &$( 1_0,20_0,11_3)$
    &$( 0_0,19_1, 6_2)$
    &$( 2_1,14_2,10_3)$\\
$Q_2$
    &$( 0_1,16_2,11_3)$
    &$( 0_0, 5_2,18_3)$
    &$( 1_0, 7_1,18_2)$
    &$( 2_0,17_1,19_3)$\\
$Q_3$
    &$( 1_0, 0_2, 3_3)$
    &$( 2_0, 1_1,10_3)$
    &$( 2_1, 5_2,20_3)$
    &$( 0_0, 3_1, 7_2)$\\
$Q_4$
    &$( 1_0, 1_1, 2_3)$
    &$( 0_1, 9_2,13_3)$
    &$( 0_0,17_1,10_2)$
    &$( 2_0,14_2,15_3)$\\
$Q_5$
    &$( 0_0,14_1,15_2)$
    &$( 1_0,13_1, 8_3)$
    &$( 2_0,10_2, 0_3)$
    &$( 0_1, 5_2,10_3)$
\end{tabular}}

{\footnotesize  \noindent
\begin{tabular}{lllllllllll}
$T'$
&$( 0_0, 1_1, 2_2, 3_3, 4_0, 5_1, 6_2)$
&$( 0_0, 3_3, 6_2, 1_1, 4_0,19_3,16_0)$
&$( 0_0, 5_1, 8_0, 3_3,18_2, 2_2,20_0)$\\
&$( 0_0,15_3, 9_1, 3_3,19_3,13_1,18_2)$
\end{tabular}}

\newpage

{\footnotesize \vspace{5pt}
\noindent {$(\alpha,\beta)=(38,3) : $}\\[5pt]
\noindent
\begin{tabular}{lllllllllll}
$P$
    &$( 3_1, 5_2,19_3)$
    &$(17_0, 5_1,12_3)$
    &$( 6_0,13_1,10_3)$
    &$( 7_1, 6_2,19_2)$
    &$(18_0, 9_1,12_2)$
    &$( 1_1, 2_1,15_2)$\\
    &$(16_0,12_1,13_3)$
    &$( 3_0, 9_2,10_2)$
    &$(18_1, 7_2,18_3)$
    &$( 5_0,13_2,17_3)$
    &$(14_0,20_1,18_2)$
    &$(12_0,13_0,10_1)$\\
    &$(10_0, 0_1, 8_3)$
    &$( 6_1,14_1, 5_3)$
    &$(15_0,20_2, 0_3)$
    &$( 1_0, 9_0, 0_2)$
    &$( 7_0,17_1, 3_2)$
    &$(11_0,15_1,11_3)$\\
    &$( 8_1, 2_2,14_3)$
    &$(16_1, 1_2, 4_3)$
    &$( 0_0,11_2,20_3)$
    &$(17_2, 2_3, 3_3)$
    &$( 4_1, 4_2, 6_3)$
    &$(20_0, 8_2, 7_3)$\\
    &$( 2_0,16_2,16_3)$
    &$( 8_0,11_1,15_3)$
    &$( 4_0,19_1,14_2)$
    &$(19_0, 1_3, 9_3)$\\
    \end{tabular}}

{\footnotesize  \noindent
\begin{tabular}{lllllllll}
$Q_1$
    &$( 0_0, 0_1, 1_2)$
    &$( 1_0,17_1,10_3)$
    &$( 2_0, 5_2, 3_3)$
    &$( 1_1,18_2,14_3)$\\
$Q_2$
    &$( 0_0, 8_1,19_2)$
    &$( 1_0, 6_1,16_3)$
    &$( 2_0,18_2,15_3)$
    &$( 1_1, 5_2,20_3)$\\
$Q_3$
    &$( 0_0,13_1,18_2)$
    &$( 1_0, 3_1,18_3)$
    &$( 2_0, 2_2, 7_3)$
    &$( 2_1,10_2, 5_3)$\\
$Q_4$
    &$( 0_0,14_1, 2_2)$
    &$( 1_0, 0_1,11_3)$
    &$( 2_0,15_2, 4_3)$
    &$( 1_1,19_2, 6_3)$\\

$T$
&$( 0_0, 7_3,11_3)$
&$( 0_0,14_2,19_3)$
&$( 0_0, 5_1,16_0)$
&$( 0_0, 2_2, 4_0)$
&$( 0_0,10_2,17_1)$
\end{tabular}}

{\footnotesize \vspace{5pt}
\noindent {$(\alpha,\beta)=(39,2) : $}\\[5pt]
\noindent
\begin{tabular}{lllllllllll}
$P$
&$(13_1,10_2,17_2)$
&$( 4_0, 7_1,20_1)$
&$( 0_0, 7_2, 0_3)$
&$(10_1,14_1, 5_3)$
&$(20_0, 4_3,18_3)$
&$(13_0, 2_1, 1_3)$ \\
&$( 6_2, 3_3,14_3)$
&$(19_0,11_1,16_2)$
&$( 1_0,17_0, 0_2)$
&$(12_1,17_1,15_3)$
&$( 3_0, 1_2, 3_2)$
&$( 3_1, 7_3,20_3)$\\
&$(11_2,15_2,10_3)$
&$(16_0, 2_3,19_3)$
&$( 2_0,15_0, 9_1)$
&$(14_0, 5_2, 8_2)$
&$(18_0, 8_1,15_1)$
&$( 6_3, 8_3,11_3)$\\
&$(11_0,20_2,12_3)$
&$( 5_0,19_2, 9_3)$
&$(12_0, 0_1,18_1)$
&$( 8_0, 9_2,14_2)$
&$( 7_0, 9_0,10_0)$
&$( 5_1,16_1, 4_2)$\\
&$( 6_0, 2_2,17_3)$
&$( 4_1, 6_1,13_3)$
&$( 1_1,18_2,16_3)$
&$(19_1,12_2,13_2)$\\
\end{tabular}}

{\footnotesize \noindent
\begin{tabular}{lllllllll}
$Q_1$
&$( 1_1, 2_1, 8_2)$
&$( 1_0, 8_0, 7_3)$
&$( 0_0, 0_1, 2_3)$
&$( 0_2,10_2, 6_3)$ \\
$Q_2$
&$( 1_1,11_3,12_3)$
&$( 0_0,14_1, 3_2)$
&$( 0_1,11_2,19_2)$
&$( 1_0, 5_0,13_3)$\\
$Q_3$
&$( 0_1, 2_2, 6_3)$
&$( 0_0,16_2,17_3)$
&$( 1_1, 9_2,19_3)$
&$( 1_0,11_0, 2_1)$\\
$Q_4$
&$( 2_0,10_1,10_2)$
&$( 1_0, 3_1,16_3)$
&$( 0_0,11_2,18_3)$
&$( 2_1,18_2, 2_3)$\\
$Q_5$
&$( 0_0,17_1,10_3)$
&$( 1_0,14_2,17_3)$
&$( 2_0, 0_1,12_2)$
&$( 1_1, 4_2, 6_3)$\\
$Q_6$
&$( 0_0, 4_1, 5_2)$
&$( 1_0, 6_1,14_3)$
&$( 2_0, 4_2,16_3)$
&$( 2_1,15_2, 3_3)$
\end{tabular}}

{\footnotesize \vspace{5pt}
\noindent {$(\alpha,\beta)=(40,1) : $}\\[5pt]
\noindent
\begin{tabular}{lllllllllll}
$P$
&$(13_1, 0_2,16_2)$
&$(10_0, 3_2, 7_3)$
&$(15_1,20_1,10_2)$
&$( 5_0, 8_0,14_2)$
&$(13_0, 2_2,17_3)$
&$( 1_0, 3_3,16_3)$ \\
&$(13_2,15_2,12_3)$
&$( 4_2, 4_3,14_3)$
&$( 3_0, 5_2,20_2)$
&$( 0_1, 3_1,15_3)$
&$( 1_1, 7_1, 6_3)$
&$(19_0,17_1,10_3)$\\
&$( 9_0, 1_2, 9_2)$
&$(19_1, 8_2,11_2)$
&$(14_1,16_1,18_3)$
&$( 0_3, 1_3,19_3)$
&$(18_0, 5_1, 9_1)$
&$( 6_0,14_0,16_0)$\\
&$( 2_0,17_0,11_1)$
&$( 4_0, 8_1,18_1)$
&$( 7_0, 6_2,20_3)$
&$(12_0, 8_3,13_3)$
&$( 0_0,18_2,19_2)$
&$(11_0, 6_1,12_2)$\\
&$( 2_1, 5_3,11_3)$
&$(20_0, 7_2, 9_3)$
&$(15_0, 4_1,12_1)$
&$(10_1,17_2, 2_3)$\\
\end{tabular}}

{\footnotesize  \noindent
\begin{tabular}{lllllllll}
$Q_1$
&$( 1_0,17_0,12_3)$
&$( 1_1, 2_1, 6_2)$
&$( 0_0, 6_1, 7_3)$
&$( 1_2, 5_2,17_3)$ \\
$Q_2$
&$( 0_1, 6_3,10_3)$
&$( 0_0,17_0, 5_3)$
&$( 1_0, 8_1, 5_2)$
&$( 1_1, 0_2,10_2)$\\
$Q_3$
&$( 1_0,12_2, 4_3)$
&$( 0_0,20_0,20_1)$
&$( 1_1, 2_2, 9_3)$
&$( 0_1,19_2,17_3)$\\
$Q_4$
&$( 1_0,12_1, 8_2)$
&$( 2_0,18_2, 0_3)$
&$( 0_0,17_1,14_3)$
&$( 1_1,13_2, 1_3)$\\
$Q_5$
&$( 2_0,14_2, 1_3)$
&$( 0_1,15_2,11_3)$
&$( 1_0, 4_1, 4_2)$
&$( 0_0, 2_1, 0_3)$\\
$Q_6$
&$( 0_0, 1_1,15_2)$
&$( 1_0,14_1, 9_3)$
&$( 2_0, 7_2, 8_3)$
&$( 0_1, 2_2, 7_3)$
\end{tabular}}

\end{document}